\documentclass[a4paper,11pt,leqno]{article}

\usepackage{amssymb}
\usepackage{graphics}
\newcommand{\color}[2][]{}

\newtheorem{lem}{Lemma}[section] \newtheorem{prp}[lem]{Proposition}
\newtheorem{df}[lem]{Definition} \newtheorem{thm}[lem]{Theorem}
\newcounter{saveitem}

\newlength{\eqdemoffset}
\newlength{\enumrqueoffset}
\newenvironment{dem}[1][] {\begin{trivlist}\item[]{\sc Proof#1.}~}
  {\ifhmode{\unskip~\nobreak}\else%
    {\nopagebreak\vspace{\eqdemoffset}\leavevmode\hfill}\fi%
    \hfill$\blacksquare$\end{trivlist}}
\newenvironment{rque}[1]{\refstepcounter{lem}%
  \begin{trivlist}\item[]{\sc #1 \arabic{section}.\arabic{lem}}~}%
  {\ifhmode{\unskip~\nobreak}\else%
    {\nopagebreak\vspace{\enumrqueoffset}\leavevmode\hfill}\fi%
    \hfill$\square$\end{trivlist}}

\renewcommand{\ss}{\scriptscriptstyle}
\newcommand{\ds}{\displaystyle}
\newcommand{\ts}{\textstyle}
\newcommand{\un}{1\!\!1}
\newcommand{\Ker}{{\rm Ker}}
\renewcommand{\Im}{{\rm Im}}
\newcommand{\Cat}{\mathcal{C}}
\newcommand{\Irr}{\mathrm{Irr}~}
\newcommand{\Mor}{\mathrm{Mor}~}
\newcommand{\Tr}{\mathrm{Tr}~}
\renewcommand{\log}{\mathrm{Log}~}
\newcommand{\Dir}{{\mathcal{D}}}
\renewcommand{\iff}{\textbf{iff} } 

\newcommand{\conv}{\makebox[.7em][c]{$\star$}}
\newcommand{\dirsum}{\raisebox{.2ex}{$\oplus$}}
\newcommand{\jok}{\bigstar}
\newcommand{\inter}[2]{[\![#1, #2]\!]}
\newcommand{\text}[1]{\mbox{#1}} 
\newcommand{\red}{\mathrm{red}}
\newcommand{\Cst}{$C^*$-\relax}
\newcommand{\tens}{\makebox[.8em][c]{$\otimes$}}
\newcommand{\tensexp}{\kern -.05em\otimes\kern -.05em}
\newcommand{\tensmax}{\tens_{\mathrm{max}}}
\newcommand{\id}{{\rm id}}
\newcommand{\rond}{\makebox[1em][c]{$\circ$}}
\newcommand{\CC}{\mathbb C} 
\newcommand{\ZZ}{\mathbb Z}
\newcommand{\NN}{\mathbb N} 
\newcommand{\FF}{\mathbb F} 
\newcommand{\Vv}{\mathfrak v}
\newcommand{\Ee}{\mathfrak e} 
\newcommand{\Gg}{\mathfrak g} 
\newcommand{\GG}{\mathfrak G} 
\newcommand{\Hh}{\mathcal{H}}
\newcommand{\Cc}{\mathfrak C} 
\newcommand{\Ss}{\mathcal{S}}
\renewcommand{\S}{E_1}
\newcommand{\T}{E_2}
\newcommand{\Te}{{\mathcal{E}}_2}
\newcommand{\Pe}{\mathcal{P}}
\renewcommand{\t}{T}
\newcommand{\Tt}{\mathcal{T}}

\title{Orientation of quantum Cayley trees \\ and applications}
\author{Roland Vergnioux} \date{}

\begin{document}

\maketitle

\begin{abstract}
  We introduce the quantum Cayley graphs associated to quantum discrete groups and study them
  in the case of trees. We focus in particular on the notion of quantum ascending orientation
  and describe the associated space of edges at infinity, which is an outcome of the
  non-involutivity of the edge-reversing operator and vanishes in the classical case. We end
  with applications to Property AO and $K$-theory.

  MSC 2000: 20G42 (58B32 46L80 19K35 46L09 46L89)
\end{abstract}

\bigskip

The original motivation of this paper is Cuntz' result on the $K$-amena\-bility of free
groups \cite{Cuntz:kmoy}, and the geometric proof of this result given by the more general
paper of Julg and Valette \cite{JulgValette:kmoy} on groups acting on trees with amenable
stabilizers.  Natural quantum analogues of the free groups are the free quantum groups
defined by Wang and van Daele \cite{DaeleWang:univ} and studied by Banica
\cite{Banic:U(n)}.  Moreover, equivariant $KK$-theory can be generalized to the case of
coactions of Hopf \Cst algebras \cite{BaajSkand:kks}, and the notion of $K$-amenability
carries over to this quantum framework without difficulty \cite{Vergnioux:these}. It is
therefore natural to ask whether free quantum groups are $K$-amenable.

To apply the method of Julg and Valette in this framework, one needs a quantum geometric object
to play the role of the tree acted upon by the quantum group under consideration.  In the case
of amalgamated free products of amenable discrete quantum groups, the construction of a quantum
analogue of the Bass-Serre tree was achieved in \cite{Vergnioux:free} and could be used to
prove the $K$-amenability of these amalgamated free products. In the case of the free quantum
groups, the needed objects should be generalizations of the Cayley graphs of the free groups.
The main goal of this paper is to define a notion of Cayley graph for discrete quantum groups,
and to study its geometric properties in the case of the free quantum groups.

We will give in the last section applications of this study: a proof of the property of Akemann
and Ostrand and the construction of a $KK$-theoretic element $\gamma$ for free quantum groups.
To prove that these quantum groups are $K$-amenable, it remains to prove that $\gamma = 1$. We
refer the reader to the last section for more historical remarks and references about Property
AO and $K$-amenability.

\bigskip

\pagebreak[3]
\noindent The paper is organized as follows: \nopagebreak
\begin{enumerate} \setlength{\itemsep}{-.5ex}
\item In the first section, we recall some notation and formulae concerning discrete quantum
  groups and classical graphs (p.~\pageref{section:notations}).
\item The second section is a technical one about fusion morphisms of free quantum groups,
  and its results are only used in the proofs of Sections~\ref{section:infinity_set}
  and~\ref{section:infinity_action}. The reader will probably like to skip over this
  section at first (p.~\pageref{section:corepr}).
\item In the third section, we give the definition of the Cayley graphs of discrete quantum
  groups and state some basic results about them (p.~\pageref{section:definition}).
\item We then restrict ourselves to the case of Cayley trees. We introduce and
  characterize this notion in the fourth section, where we also study the natural
  ascending orientation of such a tree (p.~\pageref{section:ascending}).
\item In the fifth section, we study more precisely the space of geometric edges of a quantum
  Cayley tree and we find that the projection of ascending edges onto geometric ones is not
  necessarily injective (p.~\pageref{section:geometric}).
\item We show more precisely in the sixth section that the obstruction to this injectivity
  is the existence of a natural space of (geometric) edges at infinity, which vanishes in
  the classical case (p.~\pageref{section:infinity_set}).
\item In the seventh section, we equip this space with a natural representation of the
  free quantum group under consideration, thus turning it into an interesting geometric
  object on its own (p.~\pageref{section:infinity_action}).
\item Finally the last section deals with applications, as explained above
  (p.~\pageref{section:appl}).
\end{enumerate}

\noindent {\sc Acknowledgments}

Most of the results of Sections \ref{section:corepr}--\ref{section:geometric} and, in the
special case of the orthogonal free quantum groups, of Sections \ref{section:infinity_set},
\ref{section:infinity_action} and \ref{section:gamma}, were included in my PhD thesis
\cite{Vergnioux:these} at the University Paris 7. I would like to thank my advisor, Prof.
G.~Skandalis, for having directed me to the beautiful paper \cite{JulgValette:kmoy} of Julg and
Valette, and for his precious help in the redaction of my thesis.

The generalization and continuation of these results, as well as the redaction of the present
article, were completed during a one-year stay at the University of M\"unster, where I could
benefit from a postdoctoral position. It is a pleasure to thank Prof. J.~Cuntz for his friendly
hospitality and for the very stimulating atmosphere of his team.

\section{Notation}
\label{section:notations}

The general framework of this paper will be the theory of compact quantum groups due to
Woronowicz \cite{Woro:houches}. In fact we will use it, from the dual point of view, as a
theory of discrete quantum groups.  Let us fix the notation for the rest of the paper.  The
starting object is a unital Hopf \Cst algebra $(S,\delta)$ such that $\delta(S) (1\tens S)$ and
$\delta(S) (S\tens 1)$ are dense in $S\tens S$. Such a Hopf \Cst algebra will be called a
Woronowicz \Cst algebra.  One of the key results of the theory is the existence of a unique
Haar state $h$ on $(S,\delta)$ \cite{Woro:cmp}. We will put $\delta^2 = (\id\tens\delta) \delta
= (\delta\tens\id) \delta$ and similarly $\delta^3 = (\id\tens\id\tens\delta) \delta^2$.

We denote by $\Cat$ the category of corepresentations of $(S,\delta)$ on
finite-dimensional Hilbert spaces, by $\Irr \Cat$ a set of representatives of irreducible
corepresentations modulo equivalence, and by $1_\Cat = \id_\CC\tens 1_S$ the trivial
corepresentation. We will denote by $H_\alpha$ and $v_\alpha \in B(H_\alpha)\tens S$ the
Hilbert space and the corepresentation associated to an object $\alpha\in\Cat$. The
category $\Cat$ is equipped with direct sum, tensor product and conjugation operations:
for the first and the second ones we refer to \cite{Woro:cmp}, and we give now some
precisions about the third one which is slightly more involved.

Let $(e_i)$ be an orthonormal basis of $H_\alpha$. The conjugate object $\bar\alpha$ of
$\alpha\in \Cat$ is characterized, up to isomorphism, by the existence of a conjugation
map $j_\alpha :H_\alpha \to H_{\bar\alpha}$, $\zeta \mapsto \bar\zeta$ which is an
anti-isomorphism such that ${t_\alpha : 1 \mapsto \sum e_i\tens \bar e_i}$ and $t'_\alpha :
\bar\zeta \tens \xi \mapsto (\zeta | \xi)$ are resp. elements of $\Mor (1, \alpha\tens
\bar \alpha)$ and $\Mor (\bar\alpha \tens\alpha, 1)$. We put $F_\alpha =
j_\alpha^*j_\alpha$ and we say that $j_\alpha$ is normalized if $\Tr F_\alpha = \Tr
F_\alpha^{-1}$. This positive number, which is also equal to $||t_\alpha(1)||^2$, does not
depend on the normalized map $j_\alpha$. It is called the quantum dimension of $\alpha$
and is denoted by $M_\alpha$.  When $\alpha$ is in $\Irr\Cat$, we can assume that
$\bar\alpha$ is in $\Irr\Cat$, and the possible conjugation maps $j_\alpha$ only differ by
a scalar. We have then $\bar{\bar\alpha} = \alpha$, and if $\alpha\neq \bar\alpha$ one can
choose normalized conjugation maps $j_\alpha$, $j_{\bar\alpha}$ such that $j_{\bar\alpha}
j_\alpha = 1$. If $\alpha = \bar\alpha$ one has $j_\alpha^2 = \pm 1$ for every normalized
$j_\alpha$.

The coefficients of the corepresentations $v_\alpha$ span a dense subspace $\Ss \subset S$
which turns out to be a Hopf $*$-algebra. We denote by $m$ its multiplication, and by
$\varepsilon : \Ss \to \CC$ and $\kappa : \Ss \to \Ss$ its co-unit and its antipode. Notice
that $\kappa$ is not involutive in general. In this regard, an important role is played by
a family $(f_z)_{z\in\CC}$ of multiplicative linear forms on $\Ss$, which are also related
to the non-triviality of the modular properties of $h$. We will need in this paper the
following formulae in the Hopf *-algebra $\Ss$:
\begin{eqnarray}
  \label{eq:counit}
  && \forall~ x\in\Ss~~ (\id\tens\varepsilon)\rond\delta(x) =
  (\varepsilon\tens\id)\rond\delta(x) = x \text{,} \\
  \label{eq:antipode}
  && \forall~ x\in\Ss~~ m\rond (\id\tens\kappa)\rond\delta(x) = 
  m\rond (\kappa\tens\id)\rond\delta(x) = \varepsilon(x) 1 \text{.}
\end{eqnarray}

\bigskip

Let $\Lambda_h : S \to H$ be the GNS construction of the Haar state $h$, denote by
$\lambda : S\to B(H)$ the corresponding GNS representation and by $S_\red$ its image. The
Kac system of the compact quantum group $(S,\delta)$ is given by the following formulae,
where $f\conv x:= (id\tens f)\delta(x)$ is the convolution product of $f\in \Ss^*$ and
$x\in\Ss$:
\begin{eqnarray}
  \label{eq:multunit}
  && V : \Lambda_h\tens\Lambda_h (x\tens y) \mapsto
  (\Lambda_h\tens\Lambda_h) (\delta(x) 1\tens y) \text{,} \\
  && U : \Lambda_h(x) \mapsto \Lambda_h(f_1\conv\kappa(x)) \text{.}
  \label{eq:invunit}
\end{eqnarray}
Let us recall the following notation and formulae from the general theory of multiplicative
unitaries \cite{BaajSkand:unit}. The unitary $V \in B(H\tens H)$ is multiplicative, meaning
that $V_{12}V_{13}V_{23} = V_{23}V_{12}$, and for any $\omega\in B(H)_*$ one puts $L(\omega) =
(\omega\tens\id) (V)$ and $\rho(\omega) = (\id\tens\omega) (V)$. On the other hand, $U$ is an
involutive unitary on $H$ such that $\tilde V = \Sigma (1\tens U)V(1\tens U)\Sigma$ and $\hat V
= \Sigma (U\tens 1)V(U\tens 1)\Sigma$ are again multiplicative unitaries. Moreover the
irreducibility property holds: $(\Sigma(1\tens U)V)^3 = 1$ or, equivalently, $\hat VV\tilde V =
(U\tens 1)\Sigma$. Here $\Sigma$ denotes the flip operator, and we use the leg numbering
notation.

The reduced \Cst algebra $S_\red$ coincides with the closure of $L(B(H)_*)$ in $B(H)$, and
we similarly denote by $\hat S$ the closure of $\rho(B(H)_*)$. Both can be made Hopf \Cst
algebras by the following formulae:
\begin{eqnarray}
  \label{eq:redcoprod}
  && \delta_\red (s) = V (s\tens 1) V^* = \hat V^* (1\tens s) \hat V 
  ~~\text{and} \\ \label{eq:dualcoprod}
  && \hat\delta (\hat s) = V^* (1\tens\hat s) V = \tilde V (\hat s\tens 1)
  \tilde V^* 
  \text{.} 
\end{eqnarray}
Notice that the reduction homomorphism $\lambda : S\to S_\red$ induces then an isomorphism
between the dense Hopf *-algebras of both Woronowicz \Cst algebras. Besides, the unitary $V$
lies in $M(\hat S\tens S_\red)$ and we have the following commutation relations inside $B(H)$:
$[S_\red, US_\red U] = [\hat S, U\hat SU] = 0$. There is also a full version of $S$
\cite{BaajSkand:unit} and we will say that $S$ is a full Woronowicz \Cst algebra when it
coincides with its full version.

Finally, the structure of the dual \Cst algebra $\hat S$ is very easy to describe: it is
isomorphic to the direct sum over $\alpha\in\Irr\Cat$ of the matrix \Cst algebras
$B(H_\alpha)$. We will denote by $p_\alpha \in B(H)$ the corresponding minimal central
projections of $\hat S$, except the one associated to the trivial corepresentation
$1_\Cat$ which will be denoted by $p_0$.

\bigskip

Let us recall some facts about free quantum groups. The definition was given in
\cite{Wang:freeprod,DaeleWang:univ} : let $n\geq 2$ be an integer, and $Q$ an invertible
matrix in $M_n(\CC)$, the \Cst algebra $A_u(Q)$ is then the universal unital \Cst algebra
generated by $n^2$ elements $u_{i,j}$ and the relations that make $U = (u_{i,j})$ and
$Q\bar UQ^{-1} = Q(u_{i,j}^*)Q^{-1} \in M_n(A_u(Q))$ unitary.  The \Cst algebra $A_o(Q)$
is defined similarly with the relations making $U$ unitary and $Q\bar UQ^{-1}$ equal to
$U$. We will write $S = A_o(Q)$ or $A_u(Q)$ when there is no need to distinguish the
unitary and orthogonal versions. It is easy to see that $S$ carries a unique Woronowicz
\Cst algebra structure $(S,\delta)$ for which $U$ is a corepresentation.

The corepresentation theory of $A_u(Q)$ was fully described in \cite{Banic:U(n)} in the
following way. The set of representatives $\Irr\Cat$ can be identified with the free
monoid on two generators $u$ and $\bar u$ in such a way that the corepresentation
associated to $u$ is equivalent to $U$ and the following recursive rules hold:
\begin{eqnarray*}
  &\alpha u \,\tens\, \bar u\alpha' =
  \alpha u\bar u\alpha' \oplus \alpha\tens\alpha' \text{,~~~} 
  \alpha\bar u \,\tens\, u\alpha' =
  \alpha\bar u u\alpha' \oplus \alpha\tens\alpha' \text{,} \\ 
  &\alpha u \,\tens\, u\alpha' = \alpha u u\alpha'
  \text{,~~~} \alpha\bar u \,\tens\, \bar u\alpha'=
  \alpha\bar u\bar u\alpha' \text{,~~~} 
  \overline {\alpha u} = \bar u\overline{\alpha}
  \text{,~~~} \overline {\alpha\bar u} = u\overline{\alpha} \text{.}
\end{eqnarray*}
The corepresentation theory of $A_o(Q)$ is even simpler. We assume in this case that
$Q\bar Q$ is a scalar matrix, otherwise the fundamental corepresentation $U$ is not
irreducible. The set $\Irr\Cat$ can then be identified with $\NN$ in such a
way that the corepresentation associated to $\alpha_1$ is equivalent to $U$ and the fusion
and conjugation rules read as in the representation theory of $SU(2)$:
\begin{eqnarray*}
    \alpha_k\,\tens\,\alpha_l = \alpha_{|k-l|} \oplus \alpha_{|k-l|+2} 
    \oplus \cdots \oplus \alpha_{k+l-2} \oplus \alpha_{k+l}
    \text{,~~~} \overline{\alpha_k} = \alpha_k \text{.}
\end{eqnarray*}

\bigskip

Let us finally fix some terminology concerning classical graphs.  Following
\cite{Serre:arbres}, a graph $\Gg$ will be given by a set of vertices $\Vv$, a set of edges
$\Ee$, an endpoints map $e : \Ee \to \Vv\times\Vv$ and a reversing map $\theta : \Ee\to\Ee$
which should be a involution such that $e\rond\theta = \sigma\rond e$. In this paper we denote
by $\sigma$ the flip map for spaces and \Cst algebras. If $e$ is injective, the graph $\Gg =
(\Vv,\Ee,e,\theta)$ is isomorphic to the graph $(\Vv, e(\Ee), i_{\mathrm{can}}, \sigma)$, which
we will call the simplicial realization of $\Gg$ --- although it only comes from a simplicial
complex when it has no loops, ie when $e(\Ee)$ doesn't meet the diagonal.

The set of geometric, or non-oriented, edges of $\Gg$ is the quotient $\Ee_g$ of $\Ee$ by the
relation $a \sim \theta(a)$. An orientation of the graph is a subset $\Ee_+ \subset \Ee$ such
that $\Ee$ is the disjoint union of $\Ee_+$ and $\theta(\Ee_+)$. The quotient map evidently
induces a bijection between any orientation and the set of geometric edges.  When $\Gg$ is a
tree endowed with an origin $\alpha_0$, we denote by $|\cdot|$ the distance to $\alpha_0$ and
the ascending orientation of $\Gg$ is the set of edges $a$ such that $e(a) = (\alpha,\beta)$
with $|\beta| > |\alpha|$.

Let $\Delta$ be a finite subset of a discrete group $\Gamma$ such that $1 \notin \Delta$ and
$\Delta^{-1} = \Delta$. The directional picture of the Cayley graph associated to $(\Gamma,
\Delta)$ is given by $\Vv = \Gamma$, $\Ee = \Gamma\times\Delta$, $e (\alpha,\gamma) = (\alpha,
\alpha\gamma)$ and $\theta (\alpha, \gamma) = (\alpha\gamma, \gamma^{-1})$. Its simplicial
realization will be called the simplicial picture of the Cayley graph.

\section{Complements on fusion morphisms}
\label{section:corepr}

In \cite{Banic:O(n)_cras,Banic:U(n)} a full description of the involutive semi-ring
structure of the corepresentation theory of $A_u(Q)$ and $A_o(Q)$ was given by means of
the fusion and conjugation rules on the set of irreducible objects up to equivalence. In
this section we choose concrete representatives for the irreducible objects and compute
explicitly isometric morphisms realizing the ``basic'' fusion rules. This section is a
technical one and its results are only used in Sections~\ref{section:infinity_set}
and~\ref{section:infinity_action}: we advise the reader interested in quantum Cayley
graphs to skip to the next section.

In the case of $A_u(Q)$, with $Q\in GL_n(\CC)$ and $n\geq 2$, let us choose $\gamma = u$ or
$\bar u$, and put $\gamma_{2l} = \bar\gamma$, $\gamma_{2l+1} = \gamma$. We will mainly be
interested in the corepresentations $\alpha_k = \gamma \bar\gamma \gamma \cdots \gamma_k$ ($k$
terms) and $\alpha_{k,k'} = \gamma \cdots \gamma_k \tens \gamma_{k+1}\cdots\gamma_{k+k'}$. As a
matter of fact, the fusion rules of $A_u(Q)$ reduce to the relations $\alpha_{k+1,k'+1} =
\alpha_{k+k'+2} \oplus \alpha_{k,k'}$ and trivial tensor products. In the orthogonal case, we
will also put $\gamma_k = \gamma = \alpha_1$ for every $k\in\NN$ and $\alpha_{k,k'} = \alpha_k
\tens \alpha_{k'}$, to simplify the exposition.

\bigskip

Let us now choose concrete corepresentation spaces $H_k$ and $\bar H_k$ for the classes
$\alpha_k$, $\bar\alpha_k$. We first take $H_0 = \CC$, equipped with the corepresentation
$1_\CC\tens 1_S$, and $H_\gamma = H_{\bar\gamma} = \CC^n$, equipped with the
corepresentations $U$ or $\bar U$. For any $k\in\NN^*$ we denote by $H_\gamma^{\tensexp
  k}$ the tensor product corepresentation $H_\gamma\tens H_{\bar\gamma} \tens\cdots\tens
H_{\gamma_k}$, and we define $H_k$ to be its unique sub-corepresentation equivalent to
$\alpha_k$. We proceed in the same way inside $H_{\bar\gamma_k}^{\tensexp k}$ and
$H_\gamma^{\tensexp k} \tens H_{\bar\gamma_k}^{\tensexp k'}$ to get corepresentation
spaces $\bar H_k$ and $H_{k,k'}$ representing $\bar\alpha_k$ and $\alpha_{k,k'}$. We will
denote by $t_k$, $\bar t_k$ and $t_\delta$ the morphisms associated to normalized
conjugation maps of $\alpha_k$, $\bar\alpha_k$ and $\delta\in \{\gamma, \bar\gamma\}$
respectively --- we can and will assume in this section that $j_{\bar\gamma}j_\gamma = \pm
1$, and we denote by $\mp 1$ the opposite sign. We put $m_k = M_{\alpha_k} =
M_{\bar\alpha_k}$ and we call $(m_k)_k$ the sequence of quantum dimensions of the quantum
group. Let us gather simple facts about them in the following lemma:

\begin{lem} \label{lem:conjug_morph} ~
  \begin{enumerate}
  \item \label{item:conjug_high_weight} Denote by $\t_{l} : H_\gamma^{\tensexp k} \to
    H_\gamma^{\tensexp k+2}$ the morphism $\id^{\tensexp l}\tens t_{\gamma_{l+1}}\tens
    \id^{\tensexp k-l}$. We have then $H_{k+2} = \bigcap_{l=0}^{k} \Ker~ \t^*_l \subset
    H_\gamma^{\tensexp k+2}$.
  \item \label{item:conjug_interact} We have $(\id\tens\bar t_\delta^*)(t_\delta\tens\id)
    = \pm\id_{H_\delta}$ and $t_\delta^*t_\delta = m_1 \id_\CC$ for $\delta \in \{\gamma,
    \bar\gamma\}$.
  \item \label{item:conjug_min_dim} For any $k\in\NN$ we have $m_k \geq \dim H_k$, with
    equality \iff $F_k$ is the identity. Moreover the equality $m_1 = 2$ happens only in
    the three cases $A_o\big({0\atop-1}{1\atop0}\big)$, $A_o\big({0\atop 1}{1\atop0}\big)$
    and $A_u\big({1\atop0}{0\atop1}\big)$, up to isomorphism.
  \item \label{item:conjug_induc_dim} Put $m_{-1} = 0$. The sequence of quantum dimensions
    satisfies the induction equations $m_1 m_k = m_{k+1} + m_{k-1}$ for $k\in\NN$.
    Moreover $m_0 = 1$ and $m_1$ is the geometric mean of $\Tr Q^*Q$ and $\Tr
    (Q^*Q)^{-1}$.
  \end{enumerate}
\end{lem}

\begin{dem}
  The equality of Point~\ref{item:conjug_high_weight} is true when $k=2$ because
  $H_\gamma\tens H_{\bar\gamma}$ is the orthogonal direct sum of $t_\gamma(\CC)$ and a
  subspace equivalent to $\alpha_2$, and the general result follows by induction because
  $H_{k+1} = H_{1,k} \cap H_{k,1}$. The proof of Point~\ref{item:conjug_interact} is an
  easy calculation. For Point~\ref{item:conjug_min_dim}, denote by $a$ (resp. $h$) the
  arithmetic (res. harmonic) mean of the eigenvalues of $F_k$: the normalization condition 
  of $j_k$ shows that $a=h^{-1}$ so that
  \begin{displaymath}
    m_k = a \dim H_k = \sqrt{a/h}~ \dim H_k \geq \dim H_k \text{.}
  \end{displaymath}
  For the equality case $m_1 = 2$, see \cite{Banic:U(n)}. The induction equation of
  Point~\ref{item:conjug_induc_dim} relies on the fusion rule $\alpha_k\tens\alpha_1 =
  \alpha_{k-1} \oplus \alpha_{k+1}$ which implies that $\Sigma (j_k\tens j_1)$ and
  $j_{k-1} \oplus j_{k+1}$ are normalized conjugation maps for the same corepresentation
  \cite{Woro:dual}. Finally, the formula for $m_1$ holds because the matrix $Q$ defines in
  the canonical base of $\CC^n$ a (non-normalized) conjugation map for $H_1$, by
  definition of $A_u(Q)$ \cite{Banic:U(n)}.
\end{dem}

We want now to give the explicit expression of an isometric morphism from $H_{p,p'}$
to $H_{p+1,p'+1}$, for any $p$, $p' \in \NN$. Note
that there is an evident morphism $\Tt : H_{p,p'} \to H_{p+1, p'+1}$ given by the formula
\begin{displaymath}
  \Tt = (\pi_{p+1}\tens \pi_{p'+1}) \rond (\id_{H_p} \tens t_{\gamma_{p+1}} 
  \tens \id_{H_p'}) \text{,}
\end{displaymath}
where $\pi_k$ denotes the orthogonal projection of $H_\gamma^{\tensexp k}$ onto
$H_k$. However $\Tt$ is not isometric, and its definition does not allow to compute easily 
the image of a vector $x \in H_{p,p'}$. In Proposition~\ref{prp:coeffs} we give an
explicit and simple expression of $\Tt$, which allows us to compute its polar
decomposition in Proposition~\ref{prp:polar}. From this we finally deduce 
Lemmas~\ref{lem:angles_1} and~\ref{lem:angles_2} which will be used in the proofs of
Lemmas~\ref{lem:lecture_theta} and~\ref{lem:reg_asc_commut} respectively.

\bigskip

We use more precisely the following ``basic morphisms'' from $H_\gamma^{\tensexp
  p+p'}$ to $H_\gamma^{\tensexp p+p'+2}$, indexed by $l \in \inter{0}{p}$ and $l' \in
\inter{0}{p'}$:
\begin{eqnarray*}
  && \t_{l,l'} = (\id^{\tensexp p+1}\tens t_{\gamma_{p+2}}^*\tens id^{\tensexp p'+1}) \rond
  \\ && \makebox[3cm]{}
  \rond (\id^{\tensexp p-l}\tens t_{\gamma_{p-l+1}}\tens \id^{\tensexp l+l'}\tens 
  t_{\gamma_{p+l'+1}}\tens \id^{\tensexp p'-l'}) \text{.}
\end{eqnarray*}
If $A = (a_{l,l'})$ is a $p \times p'$ matrix, we will write $\t_A = \sum a_{l,l'}
\t_{l,l'}$. Besides, we have by Lemma~\ref{lem:conjug_morph} a simpler expression of
$T_{l,l'}$ when $l$ or $l'$ equals zero:
\begin{eqnarray*}
  \t_{l,0} &=& \pm (\id^{\tensexp p-l}\tens t_{\gamma_{p-l+1}}\tens \id^{\tensexp p'+l})
  \text{~~~and} \\ \t_{0,l'} &=& 
  \pm (\id^{\tensexp p+l'}\tens t_{\gamma_{p+l'+1}}\tens \id^{\tensexp p'-l'}) \text{.}
\end{eqnarray*}

\begin{prp} \label{prp:coeffs} ~
  \begin{enumerate}
  \item \label{item:coeffs_uniq} There is at most one matrix $A$, up to a scalar factor,
    such that $\t_A$ restricts to a non-zero morphism from $H_{p,p'}$ to $H_{p+1, p'+1}$.
    If this is the case, one can assume that $a_{0,0} = 1$ and one has then $\t_A = \Tt$.
  \item \label{item:coeffs_exist} The following matrix $A$ satisfies the conditions of
    Point~\ref{item:coeffs_uniq} :
    \begin{displaymath}
      a_{l,l'} = (\mp1)^{l+l'}\, \frac{m_{p-l}m_{p'-l'}}{m_p m_{p'}} \text{.} 
    \end{displaymath}
  \end{enumerate}
\end{prp}

\begin{dem}
  \ref{item:coeffs_uniq}. It is not hard to check that the family $(\t_{l,l'})$ is free, even
  when restricted to $H_{p,p'}$. Hence it suffices to prove that an admissible $\t_A$ is
  necessarily a multiple of $\Tt$.  First of all, Point~\ref{item:conjug_high_weight} of
  Lemma~\ref{lem:conjug_morph} shows that we have $(y | \t_{l,l'}(x)) = 0$ for any $y \in
  H_{p+1,p'+1}$ and $(l,l') \neq (0,0)$. Hence if $\t_A(x) \in H_{p+1, p'+1}$ we obtain
  \begin{equation} \label{eq:simplif_norm}
    ||\t_A(x)||^2 = a_{0,0}~ (\t_A(x) | \t_{0,0}(x)) = a_{0,0}~ (\t_A(x) | \Tt(x))\text{.}
  \end{equation}
  In particular $a_{0,0}$ must be non-zero, and therefore we can assume that it equals
  $1$. To conclude we observe that the irreducible subspaces of $H_{p,p'}$ (resp. $H_{p+1,
    p'+1}$) are pairwise inequivalent, so that the morphisms $\t_A$ and $\Tt$ must be
  proportional on each irreducible subspace of $H_{p,p'}$, and (\ref{eq:simplif_norm})
  finally shows that the corresponding proportionality coefficients all equal $1$.
  
  \ref{item:coeffs_exist}. We will express the condition that $\t_A(x)$ should be in
  $H_{p+1, p'+1}$ for any $x \in H_{p,p'}$ using Point~\ref{item:conjug_high_weight} of
  Lemma~\ref{lem:conjug_morph}: for any $k \in \inter{1}{p}$ and $k' \in \inter{1}{p'}$ we
  should have $\t_{k,0}^*\t_A(x) = \t_{0,k'}^*\t_A(x) = 0$. We therefore compute, for
  $l\in \inter{0}{p}$ and $l'\in \inter{0}{p'}$:
  \begin{eqnarray*}
    \pm \t_{k,0}^*\t_{l,l'}(x) &=& (\id^{\tensexp p-k} \tens t_{\gamma_{p-k+1}}^* \tens 
    \id^{\tensexp k-1} \tens t_{\gamma_{p}}^* \tens \id^{\tensexp p'+1}) \rond \\
    && \makebox[0.5cm]{} \rond (\id^{\tensexp p-l} \tens t_{\gamma_{p-l+1}} \tens 
    \id^{\tensexp l+l'} \tens t_{\gamma_{p+l'+1}} \tens \id^{\tensexp p'-l'}) (x) \\
    &=& \left\{ 
      \begin{array}{ll}
        0 & \text{~ if $k\leq l-2$ or $k\geq l+2$} \\ 
        \pm \t_{1,0}^*\t_{0,l'}(x) & \text{~ if $k = l-1$ or $l+1$} \\ 
        m_1 \t_{1,0}^*\t_{0,l'}(x) & \text{~ if $k=l$ (see Lemma~\ref{lem:conjug_morph}).}
      \end{array}
    \right.
  \end{eqnarray*}
  The family $(\t_{1,0}^*\t_{0,l'})$ being free, we get the following conditions on $A$:
  \begin{displaymath}
    \forall~ k\in\inter{1}{p},~ l'\in\inter{0}{p'}~~
    m_1 a_{k,l'} \pm (a_{k-1,l'} + a_{k+1,l'}) = 0 \text{,}
  \end{displaymath}
  if one agrees to put $a_{p+1,l'} = 0$. We recognize the induction equations satisfied by
  the sequence $(m_{p-i}) _{0\leq i\leq p+1}$, up to a sign change. Therefore these
  conditions mean that the columns of $A$ should be proportional to $((\mp 1)^p m_p,
  \ldots, \mp m_1, 1)$.  Symmetrically the conditions $\t_{0,k'}^*\t_A(x) = 0$ are
  equivalent to the lines of $A$ being proportional to $((\mp 1)^{p'} m_{p'}, \ldots, \mp
  m_1, 1)$. The matrix of the statement satisfies these conditions, hence the associated
  morphism $\t_A$ maps $H_{p,p'}$ to $H_{p+1, p'+1}$.
\end{dem}

\begin{prp} \label{prp:polar} Let $q\in \inter{0}{\min(p,p')}$ and denote by $G\subset
  H_{p,p'}$ the subspace equivalent to $\alpha_{p+p'-2q}$. One has then
  \begin{displaymath}
    ||\Tt_{|G}||^2 = \frac{m_{p+1}m_{p'} - m_{p-q}m_{p'-q-1}} {m_p m_{p'}} \text{.}
  \end{displaymath}
\end{prp}

\begin{dem}
  Let $z \in G$ be a unit vector. Because $G$ is irreducible and $\Tt$ is a morphism, it
  is enough to compute the number $N_{p,p'}^q = ||\Tt(z)||^2$. Of course we will use the
  expression $\Tt = \t_A$ of Proposition~\ref{prp:coeffs}. We start from the
  formula~(\ref{eq:simplif_norm}) and notice that $\t_{l,l'}(z)$ is orthogonal to
  $\t_{0,0}(z) \in H_p\tens H_{\bar \gamma_p}\tens H_{\gamma_p}\tens H_{p'}$ whenever $l
  \geq 1$ or $l'\geq 1$. Hence
  \begin{eqnarray*}
    ||\t_A(z)||^2 &=& (\t_A(z) | \t_{0,0}(z)) = 
    \sum_{l, l' = 0, 1} a_{l,l'}~  (\t_{l,l'}(z) | \t_{0,0}(z)) \text{.}
  \end{eqnarray*}
  When $l$ or $l'$ equals zero, we can use the formulae for $\t_{l,l'}^*\t_{0,0}(z)$
  obtained in the proof of Proposition~\ref{prp:coeffs}. The term $l=l'=1$ will be a
  recursive one. Let us denote by $\t'_{k,k'}$ and $\Tt' = (\pi_p\tens \pi_{p'}) \rond
  \t'_{0,0}$ the morphisms analogous to $\t_{k,k'}$ and $\Tt$ for the inclusion $H_{p-1,
    p'-1} \to H_{p,p'}$. We remark that
  \begin{displaymath}
    \t_{0,0}^*\t_{1,1} = \pm \t_{0,0}^* (\id^{\tensexp p-1}\tens t_{\gamma_p}\tens 
    t_{\gamma_p} \tens \id^{\tensexp{p'-1}}) \t^{\prime *}_{0,0}
    = \t'_{0,0} \t^{\prime *}_{0,0} \text{,}
  \end{displaymath}
  so that $(\t_{1,1}(z) | \t_{0,0}(z)) = ||\t^{\prime *}_{0,0}(z)||^2 = ||\Tt^{\prime
    *}(z)||^2$. Putting all together, we get the relation
  \begin{eqnarray*}
    && N_{p,p'}^q = m_1 - \frac{m_{p-1}}{m_p} - \frac{m_{p'-1}}{m_{p'}} 
    + \frac {m_{p-1}m_{p'-1}}{m_pm_{p'}} N_{p-1,p'-1}^{q-1} \\
    &\Longleftrightarrow& m_{p'} (m_p N_{p,p'}^q - m_{p+1}) = 
    m_{p'-1} (m_{p-1} N_{p-1,p'-1}^{q-1} - m_p) \text{.}
  \end{eqnarray*}
  Hence the left-hand side quantity is invariant under simultaneous shifts of the three
  indices $p$, $p'$ and $q$. Note that the above relation is still valid when $q=0$ if one
  puts $N_{k,k'}^{-1} = 0$ for any $k$ and $k'$: as a matter of fact, in this case $z$
  lies in $H_{p+p'}$ and in particular $\Tt^{\prime *}(z) = 0$. One can therefore shift
  $q+1$ times the indices and obtain the desired identity:
  \begin{displaymath}
    m_{p'} (m_p N_{p,p'}^q - m_{p+1}) = - m_{p'-q-1} m_{p-q} \text{.}
  \end{displaymath}
\end{dem}

\begin{lem} \label{lem:angles_1}
  Let $H_{k-1,1,k} \subset H_\gamma^{\tensexp k-1} \tens H_{\gamma_k} \tens
  H_{\bar\gamma_k}^{\tensexp k}$ be the tensor product of the respective subspaces
  equivalent to $\alpha_{k-1}$, $\gamma_k$ and $\bar\alpha_k$, with $k\in\NN^*$. Let $t
  \in \Mor (H_{k-2}, H_{k-1,1})$ be an injection and denote by
  \begin{itemize}
  \item $G_1$ the subspace of $(t\tens\id) (H_{k-2, k}) \subset H_{k-1,1,k}$ equivalent to
    $\alpha_{2l}$,
  \item $G_2$ the subspace of $H_{k-1,k+1} \subset H_{k-1,1,k}$ equivalent to
    $\alpha_{2l}$.
  \end{itemize}
  Then the norm of the projection from $G_1$ onto $G_2$ equals $\sqrt{1 - \frac
    {m_lm_{l-1}} {m_km_{k-1}}}$.
\end{lem}

\begin{dem}
  Let $A = (a_{l,l'})$ and $\Tt = \t_A$ be the matrix and the morphism of
  Proposition~\ref{prp:coeffs} in the case $(p,p') = (k-2, k)$. We denote by $A' =
  (a'_{l,l'})$ the matrix given by $a'_{l,0} = a_{l,0}$ and $a'_{l,l'} = 0$ if $l'\geq 1$,
  and we remark that we have then $\t_{A'} = (\Tt'\tens\id)$, where $\Tt'$ is the morphism
  of Proposition~\ref{prp:coeffs} for $(p,p') = (k-2,0)$. Hence if $x \in
  H_\gamma^{\tensexp 2k-2}$ is a vector in the subspace of $H_{k-2, k}$ equivalent to
  $\alpha_{2l}$, we have $\t_{A'}(x) \in G_1$ and $\t_A(x) \in G_2$. The orthogonal
  projection of $G_1$ onto $G_2$ being a morphism, it is a multiple of an isometry, so
  that its norm equals
  \begin{displaymath}
    \frac{|(\t_{A'}(x) | \t_A(x))|} {||\t_{A'}(x)||\, ||\t_A(x)||} = 
    \frac{||\t_A(x)||^2} {||\t_{A'}(x)||\, ||\t_{A}(x)||} = 
    \frac {||\Tt(x)||} {||(\Tt'\tens\id)(x)||} \text{,}
  \end{displaymath}
  because the terms $T_{l,l'}(x)$ with $l'\geq 1$ are orthogonal to $\t_A(x)$. We finally
  compute the value of the last quotient thanks to Proposition~\ref{prp:polar}, with
  $(p,p',q) = (k-2,k,k-1-l)$ and $(p,p',q) = (k-2,0,0)$:
  \begin{displaymath}
    \frac{||\Tt(x)||^2} {||(\Tt'\tens\id)(x)||^2}
    = \frac{m_{k-1} m_k - m_{l-1}m_l}{m_1m_{k-2}m_k}~ \frac{m_1 m_{k-2}}{m_{k-1}}
    = 1 - \frac{m_lm_{l-1}}{m_km_{k-1}} \text{.}
  \end{displaymath}
\end{dem}

\begin{lem} \label{lem:angles_2}
  Let $H_{1,k,k'} \subset H_\gamma \tens H_{\bar\gamma}^{\tensexp k} \tens
  H_{\gamma_k}^{\tensexp k'}$ be the tensor product of the respective subspaces equivalent to
  $\gamma$, $\bar\gamma \gamma \cdots \bar\gamma_k$ ($k$ terms) and $\gamma_k \gamma_{k+1}
  \cdots \gamma_{k+k'-1}$ ($k'$ terms), with $k$, $k'\in\NN^*$. Let $t \in \Mor (H_{k-1},
  H_1\tens H_k)$ be an injection and denote by
  \begin{itemize}
  \item $G_1$ the subspace of $(t\tens\id) (H_{k-1,k'}) \subset H_{1,k,k'}$ equivalent to
    $\alpha_{k+k'-1}$, 
  \item $G_2$ the subspace of $H_{1,k+k'} \subset H_{1,k,k'}$ equivalent to $\alpha_{k+k'-1}$.
  \end{itemize}
  Then the norm of the projection from $G_1$ to $G_2$ equals $\sqrt{1 - \frac {m_{k'-1}}
    {m_{k+k'-1} m_k}}$.
\end{lem}

\begin{dem}
  Like in the previous proof we will use the morphisms $\Tt : H_{0,k+k'-1}$ $\to
  H_{1,k+k'}$ and $\Tt' : H_{0,k-1} \to H_{1,k}$ studied in Proposition~\ref{prp:coeffs}.
  We notice that $G_1$ (resp. $G_2$) is the image of $H_{k+k'-1}$ by $(\Tt'\tens\id)$
  (resp. $\Tt$), so that the norm of the projection we are interested in is given by
  \begin{displaymath}
    \frac{|(\Tt(x) | (\Tt'\tens\id)(x))|} {||\Tt(x)||\, ||(\Tt'\tens\id)(x)||} =
    \frac{||\Tt(x)||} {||(\Tt'\tens\id)(x)||} \text{,}
  \end{displaymath}
  for the same reason as above. Proposition~\ref{prp:polar} with $(p,p',q) = (0,k+k'-1,0)$ and
  $(0,k-1,0)$ gives then
  \begin{displaymath}
    \frac{||\Tt(x)||^2} {||(\Tt'\tens\id)(x)||^2} = \frac{(m_1m_{k+k'-1} - m_{k+k'-2})m_{k-1}} 
    {m_{k+k'-1} (m_1m_{k-1} - m_{k-2})} = \frac {m_{k+k'}m_{k-1}} {m_{k+k'-1} m_k} \text{.}
  \end{displaymath}
  The result follows then from the identity $m_{k+k'-1} m_k = m_{k+k'} m_{k-1} + m_{k'-1}$,
  which is easy to prove by induction, or by noticing that the irreducible subobjects of
  $H_{k+k'-1,k}$ are the same as for $H_{k+k', k-1}$, up to the one equivalent to $H_{k'-1}$.
\end{dem}

\section{Quantum Cayley graphs}
\label{section:definition}

In this section we introduce the notion of Cayley graph for discrete quantum groups. In
fact the classical notion can be generalized into two different directions, coming from
the two different pictures introduced in Section~\ref{section:notations}. The quantum
generalization of the simplicial picture is still a classical graph. On the contrary, the
$\ell^2$-spaces of the directional picture give rise in the quantum case to a quantum
object, in the spirit of non-commutative geometry.

\bigskip

In the following definition, we use freely the notation of Section~\ref{section:notations}. In
particular, $S$ and $\hat S$ are the dual Hopf \Cst algebras of a compact quantum group --- $S$
being unital ---, $H$ is the GNS space of the Haar state of $S$ and $p_\alpha$ is the minimal
central projection of $\hat S$ corresponding to an irreducible corepresentation
$\alpha\in\Irr\Cat$.

\pagebreak[3]
\begin{df} \label{df:cayley}
  Let $S$ be a Woronowicz \Cst algebra and $p_1$ a central projection of $\hat S$ such that
  $Up_1U = p_1$ and $p_0 p_1 = 0$. \nopagebreak
  \begin{enumerate}
  \item The classical Cayley graph $\Gg$ associated with $(S,p_1)$ is given in simplicial
    form by $\Vv = \Irr\Cat$ and $\Ee = \{ (\alpha,\alpha') \in \Vv^2 ~|~
    \hat\delta(p_{\alpha'}) (p_\alpha\tens p_1) \neq 0\}$.
  \item The hilbertian quantum Cayley graph associated with $(S,p_1)$ is the $4$-uplet
    $(H, K, E, \Theta)$ where $K = H\tens p_1 H$, $E = V_{|K}\in B(K, H\tens H)$ and
    $\Theta = \tilde V(1\tens U)_{|K} \in B(K)$.  \setcounter{saveitem}{\value{enumi}}
  \end{enumerate}
  Let us introduce some more objects associated with this quantum graph. We denote by
  $\epsilon$ be the linear form on $p_1 H$ defined by $\epsilon(\Lambda_h(x)) =
  \Lambda_h(\varepsilon(x))$.
  \begin{enumerate}
    \setcounter{enumi}{\value{saveitem}}
  \item The source and target operators of the hilbertian quantum Cayley graph are $\S =
    (\id\tens\epsilon)$ and $\T = \S\rond\Theta \in B(K,H)$.
  \item The quantum $\ell^2$-space of geometric edges is $K_g = \Ker (\Theta + \id)$.
  \end{enumerate}
\end{df}

\begin{rque}{Remarks} \label{rque:cayley}
  \begin{enumerate}
  \item \label{enum:fond_set} The central projections $p_1$ that match the hypotheses of
    Definition~\ref{df:cayley} are sums of projections $p_\alpha$ over finite subsets
    $\Dir \subset \Irr\Cat$ such that $\bar\Dir = \Dir$ and $1_\Cat \notin\Dir$.  The
    elements of $\Ee$ are then the ordered pairs of vertices $(\alpha,\alpha')$ for which
    there exist $\gamma\in\Dir$ such that $\alpha' \subset \alpha\tens\gamma$. Note that
    this set of edges is symmetric, thanks to the equivalence $\alpha' \subset
    \alpha\tens\gamma \Leftrightarrow \alpha\subset \alpha'\tens\bar\gamma$ (Jacobi
    duality).  In this paper, the classical Cayley graph will mainly be used as a tool for
    the study of the quantum one.
  \item The hilbertian quantum Cayley graph will be more useful for our purposes because
    he naturally carries representations of the discrete quantum group under
    consideration: the \Cst algebra $S$ acts on $H$ via the GNS representation, and we let
    it act trivially on $p_1 H$. Moreover the operators $\Theta$, $\S$ and $\T$ commute to
    these representations, and in particular $K_g$ is also endowed with a natural
    representation of $S$. The commutation properties to the action of $\hat S$ will be
    examined in Proposition~\ref{prp:commutdual}.
  \item \label{enum:non_invol} The identity $\hat VV\tilde V = (U\tens 1) \Sigma$ provides
    us with another expression for the reversing operator: $\Theta = (\Sigma\hat VV)^*$.
    Moreover the identity $\tilde V^* = (\hat J\tens J) \tilde V (\hat J\tens J)$, where
    $J$, $\hat J$ are the modular conjugations of $S$ and $\hat S$ \cite{KustVaes:lcqg},
    shows that$(\hat J\tens\hat J) \Theta (\hat J\tens\hat J) = \Theta^* = \Theta^{-1}$.
    But the main fact about the reversing operator is its non-involutivity in the quantum
    case, see Proposition~\ref{prp:thetaclassic} --- in fact in this proposition it is
    enough to consider the restriction of $\tilde V (1\tens U)$ to $H\tens p_1H$, as soon
    as $\Dir$ generates $\Cat$.
  \end{enumerate} 
\end{rque}

\begin{rque}{Example}(classical case) \label{rque:classical}
  Suppose $S = C^*(\Gamma)$ for some discrete group $\Gamma$, with the co-commutative coproduct
  given by $\delta(\gamma) = \gamma\tens\gamma$ for all $\gamma\in\Gamma \subset C^*(\Gamma)$.
  Then $\Irr\Cat$ identifies with $\Gamma$ in such a way that $v_\gamma \simeq
  \id_\CC\tens\gamma$ for every $\gamma\in\Gamma$, and the tensor product (resp.  the
  conjugation) of corepresentations then coincides with the product (resp. the inverse) of
  $\Gamma$. In particular, the inclusion $\alpha'\subset \alpha\tens \gamma$ reduces in this
  case to an equality $\alpha' = \alpha\gamma$, so that the classical graph of
  Definition~\ref{df:cayley} is nothing but the simplicial picture of the Cayley graph
  associated to $(\Gamma, \Delta)$.
  
  Besides, one has $H = \ell^2(\Gamma)$, $S_\red = C^*_\red(\Gamma)$, $\hat S =
  c_0(\Gamma)$ and the projections $p_\alpha$ correspond to the characteristic functions
  $\un_\alpha \in \hat S$ of the points of $\Gamma$. Moreover one has the following
  expressions for the Kac system of $(S,\delta)$: $V (\un_\alpha \tens \un_\beta) =
  \un_\alpha \tens \un_{\alpha\beta}$ and $U(\un_\alpha) = \un_{\alpha^{-1}}$. From this
  it is easy to see that the hilbertian quantum graph of Definition~\ref{df:cayley} is
  nothing but the $\ell^2$-object associated to the directional picture $(\Vv, \Ee, e,
  \theta)$ of the Cayley graph of $(\Gamma, \Delta)$. The only non-trivial check concerns
  the reversing operator: according to Proposition~\ref{prp:thetaclassic}, one has
  $\Sigma\hat VV = V^*\Sigma V = E^*\Sigma E$ so that $\Theta^* = \Theta$ is the classical
  reversing operator.
\end{rque}

\begin{prp} \label{prp:thetaclassic} Let $(H,V,U)$ be an irreducible
  Kac system \cite{BaajSkand:unit}. Then the multiplicative unitary $V$ is co-commutative
  \iff $\hat V = \Sigma V^* \Sigma$ \iff $\tilde V (1\tens U)$ is involutive.
\end{prp}

\begin{dem}
  The direct implications are easy to check in the underlying locally compact groups.
  Conversely, assume that $\hat V = \Sigma V^* \Sigma$. Then for any $x = (\id\tens\omega)
  (V) \in \hat S_\red$, one also has $x = U (\id\tens\omega) (V^*) U \in U\hat S_\red U
  \subset \hat S'_\red$, hence $\hat S_\red$ is commutative. Replacing $V$ by $\tilde V$
  one gets the dual version of this result: if $\tilde V = \Sigma V^*\Sigma$, then $V$ is
  commutative. Now, $\tilde V (1\tens U)$ is involutive \iff $\tilde V (1\tens U) =
  (1\tens U) \tilde V^*$ \iff $\Sigma (1\tens U)\tilde V (1\tens U)\Sigma = \Sigma \tilde
  V^* \Sigma$, which implies by the previous ``dual'' statement that $\tilde V$ is
  commutative, hence $V$ is co-commutative.
\end{dem}

Let us give now alternative expressions for the reversing, source and target operators in
terms of the Hopf $*$-algebra structure of $\Ss$, and study the intertwining properties of
these operators relatively to the representations of the dual Hopf \Cst algebra $\hat S$.

\begin{lem} \label{lem:algTheta} Let $x$, $y\in \Ss \subset S_\red$,
  we have $\tilde V (1\tens U)\rond (\Lambda_h\tens\Lambda_h) (x\tens y) =$ 
  $(\Lambda_h\tens\Lambda_h) \rond (\id\tens\kappa) ((x\tens 1) \delta(y))$.
\end{lem}

\begin{dem}
  In this proof we will write $\Theta$ in place of $\tilde V(1\tens U)$, although we do
  not necessarily restrict ourselves to $K$. We have $\Theta \rond
  (\Lambda_h\tens\Lambda_h) (x\tens y) = (x\tens \id) \rond \Theta \rond
  (\Lambda_h\tens\Lambda_h) (1\tens y)$, so that it suffices to consider the case when
  $x=1$. Let us use the expressions~(\ref{eq:multunit}) and~(\ref{eq:invunit}) of $U$ and
  $V$:
  \begin{eqnarray*}
    \Theta \rond(\Lambda_h\tens\Lambda_h) (1\tens y) &=&
    \Sigma(1\tens U)V(1\tens U) \rond(\Lambda_h\tens\Lambda_h) 
    (f_1\conv\kappa(y)\tens 1) \\
    &=& \Sigma(1\tens U) \rond(\Lambda_h\tens\Lambda_h) \rond
    \delta(f_1\conv\kappa(y)) \text{.}
  \end{eqnarray*}
  It is easy to check that $\delta(f_z\conv a) = (\id\tens(f_z\conv)) (\delta(a))$, hence
  \begin{eqnarray*}
    \Theta \rond(\Lambda_h\tens\Lambda_h) (1\tens y) &=&
    \Sigma(1\tens U) \rond(\Lambda_h\tens\Lambda_h)
    (\id\tens (f_1\conv)) (\delta(\kappa(y)))  \\
    &=& \Sigma(1\tens U) \rond(\Lambda_h\tens\Lambda_h)
    (\kappa\tens (f_1\conv)\kappa)\sigma\delta(y)\text{.}
  \end{eqnarray*}
  One recognizes then $(1\tens U)^2 = 1\tens 1$:
  \begin{eqnarray*}
    \Theta \rond(\Lambda_h\tens\Lambda_h) (x\tens y) &=& 
    \Sigma \rond(\Lambda_h\tens\Lambda_h) \rond 
    (\kappa\tens\id) \sigma\delta(y) \\ 
    &=& (\Lambda_h\tens\Lambda_h) \rond 
    (\id\tens\kappa) \delta(y) \text{.} 
  \end{eqnarray*}
\end{dem}

\begin{prp} \label{prp:endpoints}
  Let $\epsilon$ be the linear form on $\Lambda_h(\Ss)$ defined by $\epsilon\rond\Lambda_h
  = \Lambda_h\rond\varepsilon$.  We have the following identities:
  \begin{enumerate}
  \item $\S = (\id\tens\epsilon)\rond V$ and $\T = (\epsilon\tens\id)\rond V$ on $K\cap
    (\Lambda_h\tens\Lambda_h) (\Ss\tens\Ss)$,
  \item $\T\rond (\Lambda_h\tens\Lambda_h) (x\tens y) = \Lambda_h(xy)$ for $x\tens y\in
    \Ss\tens\Ss$, and $\T\rond \Theta = \S$.
  \end{enumerate}
\end{prp}

\begin{dem}
  The co-unit $\varepsilon$ being multiplicative, we have for all $x$ and $y$ in $\Ss$
  \begin{displaymath}
    (\id\tens\varepsilon)(\delta(x)(1\tens y)) = 
    \varepsilon(y) (\id\tens\varepsilon) \delta(x) = 
    \varepsilon(y) x = (\id\tens\varepsilon)(x\tens y) \text{,}
  \end{displaymath} 
  hence $\S = (\id\tens\epsilon)\rond V$. In the same way one can write, using the
  identity $\varepsilon\rond\kappa = \varepsilon$ and Equation~(\ref{eq:counit}):
  \begin{eqnarray*}
    (\id\tens\varepsilon)(\id\tens\kappa)((x\tens 1)\delta(y)) &=& 
    (\id\tens\varepsilon)((x\tens 1)\delta(y)) \\ &=& xy ~=~ 
    (\varepsilon\tens\id)(\delta(x)(1\tens y)) \text{,}
  \end{eqnarray*}
  which yields $\T\rond (\Lambda_h\tens\Lambda_h) (x\tens y) = \Lambda_h(xy)$ and $\T =$
  $(\epsilon\tens\id)\rond V$, thanks to the definition of $\T$ and the expression of
  $\Theta$ given by Lemma~\ref{lem:algTheta}.  Now, using these results and
  Equation~(\ref{eq:antipode}), we can proceed to the last computation, where $m :
  \Ss\tens\Ss \to \Ss$ denotes the multiplication of $\Ss$:
  \begin{eqnarray*}
    \T\rond\Theta\rond (\Lambda_h\tens\Lambda_h) (x\tens y) &=&
    \Lambda_h( x~ (m(\id\tens\kappa)\delta(y)) ) \\
    &=& \varepsilon(y)~ \Lambda_h(x) ~=~
    \S\rond (\Lambda_h\tens\Lambda_h) (x\tens y) \text{.}
  \end{eqnarray*}
\end{dem}

\begin{prp} \label{prp:commutdual}
  Let us define $\hat\pi_2 : \hat S^{\tensexp 2} \to L(H)$ by the formula
  $\hat\pi_2(x\tens x') = x (Ux'U)$. Similarly, let us denote by $\hat\pi_4 : \hat
  S^{\tensexp 4} \to L(K)$ the homomorphism such that $\hat\pi_4(x\tens y\tens y'\tens x')
  = (x\tens y) (Ux'U\tens Uy'U)$, and let us put $\hat\delta' = \hat\pi_4 \rond (1\tens
  1\tens\hat\delta)$, so that $\hat\delta'(x) = (U\tens U)\Sigma
  \hat\delta(x)\Sigma(U\tens U)$. One has then, for any $x\in \hat S$: \pagebreak[3]
  \begin{enumerate}
  \item \label{enum:commutdual1} $\Theta \rond (x\tens 1) = \hat\delta(x) \rond \Theta$,
  \item \label{enum:commutdual2} $\Theta \rond (1\tens x) = (1\tens UxU) \rond \Theta$,
  \item \label{enum:commutdual3} $\Theta \rond \hat\delta'(x) = (UxU\tens 1) \rond \Theta$,
  \item \label{enum:commutdualT} $\T\rond \hat\delta(x) = x \rond\T$ and $\T\rond
    \hat\delta'(x) = UxU \rond\T$.
  \end{enumerate}
  Hence $\Theta$ intertwines the representations $\hat\pi_4\rond (\id\tens\id\tens\hat\delta)$
  and $\hat\pi_4\rond (\hat\delta\tens\id\tens\id)$ of $\hat S\tens\hat S\tens\hat S$ on $K$.
  In particular $\Theta$ commutes to $\hat\pi_4 \rond \hat\delta^3$. Similarly, $\T$ intertwines
  the representations $\hat\pi_2$ and $\hat\pi_4 \rond (\hat\delta\tens\hat\delta)$ of $\hat
  S\tens\hat S$.
\end{prp}

\begin{dem}
  Point~\ref{enum:commutdual1} results from the identity $\hat\delta(x) = \tilde V(x\tens
  1)\tilde V^*$. Writing $\T = (\id\tens\epsilon)\rond\Theta^{-1}$, it implies the first
  relation of point~\ref{enum:commutdualT}. For point~\ref{enum:commutdual3}, one uses the
  formula $\Theta = (\Sigma\hat V V)^*$ and the fact that $V$ commutes to $U\hat SU\tens 1$:
  \begin{eqnarray*}
    (UxU\tens 1)\Theta &=& (UxU\tens 1) V^* \hat V^* \Sigma
    = V^* (UxU\tens 1) \hat V^* \Sigma \\ 
    &=& V^* \hat V^* (U\tens U)\tilde V (x\tens 1)
    \tilde V^* (U\tens U)\Sigma \\ 
    &=& V^* \hat V^* (U\tens U) \hat\delta(x)(U\tens U)\Sigma
    = V^* \hat V^* \Sigma \hat\delta'(x) 
    = \Theta \hat\delta'(x) \text{.}
  \end{eqnarray*}
  Composing on the left by $(\id\tens\epsilon)$, one obtains the second relation of
  point~\ref{enum:commutdualT}. For point~\ref{enum:commutdual2}, simply notice that $\tilde V$
  is in $M(USU\tens\hat S)$, and hence commutes to $1\tens U\hat SU$.
\end{dem}

\section{Ascending orientation}
\label{section:ascending}

In the case of a classical tree, the ascending orientation associated to a chosen origin
defines a subspace $K_+$ of the $\ell^2$-space of edges $K$. The aim of this section is to
introduce and study such a subspace in the case of quantum Cayley graphs. The next
definition relies on the links between the quantum and classical Cayley graphs, the latter
one being endowed with the origin $1_\Cat$.

\begin{df} \label{df:orient}
  Let $S$ be a Woronowicz \Cst algebra and $p_1$ a central projection of $\hat S$ such that
  $Up_1U = p_1$ and $p_0 p_1 = 0$. Assume that the classical Cayley graph associated with
  $(S,p_1)$ is a tree, and denote by $|\,\cdot\,|$ the distance to the origin $1_\Cat$ in this
  tree.
  \begin{enumerate}
  \item For any $n\in \NN \setminus \{0,1\}$ we put $p_n = \sum \{p_\alpha ~|~ |\alpha| = n\}
    \in Z(\hat S)$.
  \item We call $p_{\ss\jok+} = \sum\, (p_n\tens p_1) \hat\delta(p_{n+1})$ and
    $p_{\ss+\jok} = \sum\, (p_n\tens p_1) \hat\delta'(p_{n+1})$ the left and right
    ascending projections. Put $p_{\ss\jok-} = 1 - p_{\ss\jok+}$, $p_{\ss-\jok} = 1 -
    p_{\ss+\jok}$.
  \item We call $p_{\ss++} = p_{\ss+\jok} p_{\ss\jok+}$ the ascending projection of the
    quantum Cayley tree, and we denote by $K_{\ss++}$ its image. We define similarly
    \begin{eqnarray*}
      p_{\ss+-} = p_{\ss+\jok} p_{\ss\jok -} \text{,~~~}
      &p_{\ss-+} = p_{\ss-\jok} p_{\ss\jok +} \text{~~~and~~~} 
      &p_{\ss--} = p_{\ss-\jok} p_{\ss\jok -} \text{,} \\
      K_{\ss+-} = p_{\ss+-} K \text{,~~~}
      &K_{\ss-+} = p_{\ss-+} K \text{~~~and~~~}
      &K_{\ss--} = p_{\ss--} K \text{.}
    \end{eqnarray*}
  \end{enumerate}
\end{df}

\pagebreak[3]
\begin{rque}{Remarks} \nopagebreak
  \begin{enumerate}
  \item We have $|\alpha| = 1 \Leftrightarrow \alpha \in \Dir$ and $|\alpha| = 0
    \Leftrightarrow \alpha = 1_\Cat$. In particular the first point of
    Definition~\ref{df:orient} is consistent with the notation $p_0$ and $p_1$ used in
    Definition~\ref{df:cayley}.
  \item \label{enum:rque_p++_orient} Take $\alpha\in\Irr\Cat$ with $|\alpha| = n+1$.  We
    have $\hat\delta(p_\alpha) = \sum \hat\delta(p_\alpha) (p_\beta\tens p_{\beta'})$,
    where the sum goes over the ordered pairs $(\beta, \beta')$ such that $\alpha\subset
    \beta\tens\beta'$. Hence $\hat\delta(p_\alpha) (p_n\tens p_1) = \hat\delta(p_\alpha)
    (p_{\alpha'}\tens p_1)$, where $\alpha'$ is the vertex preceding $\alpha$ in the
    classical Cayley graph $\Gg$. Hence we get the following expression of $p_{\ss\jok+}$,
    in terms of the classical ascending orientation $\Ee_+$ of $\Gg$:
    \begin{displaymath} \label{eq:p*+_geom}
      p_{\ss\jok+} = \sum_{(\alpha',\alpha) \in \Ee_+} V^* (p_{\alpha'}\tens p_\alpha) V 
      (\id\tens p_1) \text{.}
    \end{displaymath}
    Recall that $V$ plays the role of the endpoints operator, which implements in the
    co-commutative case the equivalence between the simplicial and directional pictures of
    the Cayley graph.
  \item Let $J$ (resp. $\hat J$) be the modular conjugation on $H$ induced by the
    involution of $S$ (resp. $\hat S$). We know from \cite{KustVaes:lcqg} that $U = \hat
    JJ$, $[J, p_n] = 0$ and $(J\tens J) \hat\delta(p_n) (J\tens J) =
    \Sigma\hat\delta(p_n)\Sigma$. From this we can deduce the following relation between
    $p_{\ss+\jok}$ and $p_{\ss\jok+}$:
    \begin{displaymath} \label{eq:p*+_p+*}
      p_{\ss+\jok} = (\hat J\tens\hat J) p_{\ss\jok+} (\hat J\tens\hat J) \text{.}
    \end{displaymath}
    Hence $p_{\ss+\jok}$ and $p_{\ss\jok+}$ come from the same projection of $M(\hat
    S\tens\hat S)$ acting respectively on the left and on the right of $K$. In particular
    they commute and are equal in the co-commutative case.
  \end{enumerate}
\end{rque}

In the next proposition we examine the links between the ascending projections and the
reversing and target operators. The first point of the proposition shows that the
reversing operator switches the left and right versions of the quantum ascending
projections: this is the reason why it is necessary to use both $p_{\ss\jok+}$ and
$p_{\ss+\jok}$ in the general case.
  
\begin{prp} \label{prp:orient} 
  We use the hypothesis and notation of Definition~\ref{df:orient}.
  \begin{enumerate}
  \item \label{enum:rev_orient} We have $p_{\ss\jok-} = \Theta p_{\ss+\jok} \Theta^*$ and
    $p_{\ss-\jok} = \Theta^* p_{\ss\jok+} \Theta$. More precisely:
    \begin{eqnarray*}
      &&\Theta p_{\ss+\jok} (p_n\tens \id) = 
      (p_{n+1}\tens \id) p_{\ss\jok-} \Theta \text{~~~and~~~} \\
      &&\Theta p_{\ss-\jok} (p_n\tens \id) = 
      (p_{n-1}\tens \id) p_{\ss\jok+} \Theta \text{.}
    \end{eqnarray*}
  \item \label{enum:target_orient} We have $\T p_{\ss+-} = \T p_{\ss-+} = 0$ and
    \begin{displaymath}
      p_n \T = \T (p_{n-1}\tens \id) p_{\ss++} + 
      \T (p_{n+1}\tens \id) p_{\ss--} \text{.}
    \end{displaymath}
  \end{enumerate}
\end{prp}

\begin{dem}
  \ref{enum:rev_orient}. We put $u = \mathrm{Ad}(U)$.  Using the formulae $\tilde V
  (p_\alpha\tens 1) \tilde V^* = \hat\delta(p_\alpha)$ and $V^*(1\tens p_\alpha)V =
  \hat\delta(p_\alpha)$ for the dual coproduct, and the fact that $p_n$ commutes to
  $U$, we see that
  \begin{eqnarray*}
    \tilde V^*(p_n\tens 1)\tilde V &=& (U\tens 1) \Sigma V^* \Sigma 
    (Up_n U\tens 1) \Sigma V \Sigma (U\tens 1) \\ 
    &=& (u\tens\id)\sigma (V^* (1\tens p_n) V) = 
    (u\tens\id)\sigma\hat\delta(p_n) \text{.}
  \end{eqnarray*}
  We use this expression in conjunction with the definition of $\Theta$:
  \begin{eqnarray*}
    \Theta p_{\ss+\jok}(p_n\tens 1) \Theta^* &=& \tilde V 
    (p_n\tens p_1) (u\tens\id)\sigma\hat\delta(p_{n+1}) \tilde V^* 
    \\ &=& (1\tens p_1) \hat\delta(p_n) \tilde V 
    (u\tens\id) \sigma\hat\delta(p_{n+1}) \tilde V^* 
    = \hat\delta(p_n) (p_{n+1}\tens p_1) \text{.}
  \end{eqnarray*}
  To prove that the last expression is equal to $p_{\ss\jok-} (p_{n+1}\tens p_1)$, it is
  enough to check that we obtain $p_{\ss\jok-}$ by summing it over $n$. But $(p_n\tens
  p_1) \hat\delta(p_{n'})$ vanishes as soon as $n'\neq n\pm1$, so that
  \begin{eqnarray*}
    \sum (p_{n+1}\tens p_1)\hat\delta(p_n) + p_{\ss\jok +} &=& 
    \sum \big((p_{n+1}\tens p_1) \hat\delta(p_n) + 
    (p_n\tens p_1) \hat\delta(p_{n+1})\big) \\ 
    &=& \big(\ts\sum (p_n\tens p_1)\big) 
    \big(\ts\sum \hat\delta(p_{n'})\big) = 
    \id_K \text{.}
  \end{eqnarray*}
  
  \ref{enum:target_orient}. The last point of Proposition~\ref{prp:commutdual} shows that
  $p_k \T (p_l\tens p_1)$ equals $\T \hat\delta(p_k) (p_l \tens p_1)$. In particular $\T
  (p_n\tens p_1) p_{\ss\jok +} = p_{n+1} \T p_{\ss\jok+}$ and similarly $\T (p_n\tens
  p_1)p_{\ss\jok -} = p_{n-1} \T p_{\ss\jok -}$. On the other hand one has, using
  Propositions~\ref{prp:endpoints} and \ref{prp:orient}:
  \begin{eqnarray*}
    \T (p_n\tens p_1) p_{\ss-\jok} &=&
    \S \Theta (p_n\tens p_1) p_{\ss-\jok} =
    \S (p_{n-1}\tens p_1) p_{\ss\jok +} \Theta \\ 
    &=& p_{n-1} \S p_{\ss\jok +} \Theta =
    p_{n-1} \T p_{\ss-\jok} \text{,~~~and similarly} \\
    \T (p_n\tens p_1) p_{\ss+\jok} &=& 
    p_{n+1} \T p_{\ss+\jok} \text{.}
  \end{eqnarray*}
  As a result $\T (p_n\tens p_1) p_{\ss-+}$ equals both $p_{n+1} \T p_{\ss-+}$ and
  $p_{n-1} \T p_{\ss-+}$, so that it must vanish, and in the same way $\T (p_n\tens p_1)
  p_{\ss+-} = 0$. In particular $p_n \T = p_n \T {(p_{\ss++} + p_{\ss--})}$, and the last
  statement of the proposition results then from the identities $p_n \T p_{\ss\jok+} = \T
  (p_{n-1}\tens p_1)p_{\ss\jok+}$ and $p_n \T p_{\ss\jok-} = \T (p_{n+1}\tens
  p_1)p_{\ss\jok-}$ that we proved above.
\end{dem}

Until now we have used a very minimal notion of ``tree'' for our quantum Cayley graphs,
namely the fact that the corresponding classical Cayley graph should be a classical tree.
However this notion is too weak for our purposes, because it doesn't take into account
multiplicity issues that appear in the quantum case. More precisely, let us define the
``full'' classical Cayley graph $\GG$ associated to $(S, p_1)$ in the following way:
\begin{eqnarray*}
  && \Vv = \Irr\Cat \text{,~~~} \Ee = \{ (\alpha,\alpha',\gamma,i) ~|~ \gamma\in\Dir,~ 
  \alpha'\subset \alpha\tens\gamma \textrm{~with mult. order $i$} \}  \\
  && e(\alpha,\alpha',\gamma,i) = (\alpha,\alpha') \text{,~~~} 
  \theta(\alpha,\alpha',\gamma,i) = (\alpha',\alpha,\bar\gamma,i) \text{.}
\end{eqnarray*}
Here $\Dir$ stands for the set of corepresentations associated with $p_1$, like in
Remark~\ref{rque:cayley}.\ref{enum:fond_set}. The image of $\GG$ by $e$ is the classical
Cayley graph $\Gg$ of Definition~\ref{df:cayley}, but the map $e$ needs not to be
injective in general. The component $\gamma$ of an edge $(\alpha,\alpha',\gamma,i)$ is
called the direction of the edge. In the rest of this paper, we will assume that the full
classical Cayley graph $\GG$ with origin $1_\Cat$ is a ``directional tree'', meaning that
it is a tree and that the ascending edges starting from a given vertex have pairwise
different directions.

In Lemma~\ref{lem:strict} we state some basic results about classical Cayley graphs and
give a corepresentation-theoretic formulation of the extra assumptions introduced above.
Proposition~\ref{prp:but_inj} shows that our framework is the right one for the study of
free quantum groups, ie free products of orthogonal and unitary free quantum groups
\cite{DaeleWang:univ,Banic:U(n)}. Finally we prove that the quantum ascending orientation
$K_{\ss++} \subset K$ behaves nicely in this framework: the target operator induces a
bijection between ascending edges and vertices orthogonal to the origin, exactly like in
the classical case.

\begin{lem} \label{lem:strict}
  Let $S$ be a Woronowicz \Cst algebra and $p_1$ a central projection of $\hat S$ such
  that $Up_1U = p_1$ and $p_0 p_1 = 0$. Assume that the classical Cayley graph $\Gg$ is a
  tree and denote by $(\alpha\tens \gamma)_+$ (resp.  $(\alpha\tens\gamma)_-$) the sum of
  the subobjects of $(\alpha\tens \gamma)$ which are further from (resp. closer to)
  $1_\Cat$ than $\alpha$.
  \begin{enumerate}
  \item \label{item:sym_dist} For every $\alpha\in\Irr\Cat$ one has $|\alpha| =
    |\bar\alpha|$ in $\Gg$.
  \item \label{item:strict_repr} The full classical Cayley graph $\GG$ is a directional tree
    \iff
    \begin{itemize}
    \item for all $\alpha\in\Irr\Cat$ and $\gamma\in\Dir$, $(\alpha\tens\gamma)_+$ is
      irreducible or zero and
    \item for all $\alpha\in\Irr\Cat$ and $\gamma\neq\gamma' \in\Dir$, $(\alpha\tens\gamma)_+$
      and $(\alpha\tens\gamma')_+$ are inequivalent or zero.
    \end{itemize}
  \item \label{item:always_asc} We assume that $\GG$ is a directional tree. For any
    $(\alpha,\gamma) \in \Irr\Cat\times\Dir$, one has $(\alpha\tens \gamma)_+ = 0$ \iff
    $\dim\gamma = 1$ and $\alpha$ is the target of an ascending edge with direction
    $\bar\gamma$.
  \item \label{item:dim_grows} If $\GG$ is a directional tree and $(\alpha,\beta)$ is an ascending
    edge then $\dim\beta \geq \dim\alpha$, with equality \iff the corresponding direction
    $\gamma\in\Dir$ has dimension $1$.
  \end{enumerate}
\end{lem}

\begin{dem}
  \ref{item:sym_dist}. For this first point $\Gg$ does not need to be a tree. Because
  $\bar{\Dir} = \Dir$, it is enough to prove the following property: $|\alpha| \leq n$ \iff
  there exist elements $\gamma_1$, \ldots, $\gamma_n \in \Dir$ such that $\alpha \subset
  \gamma_1\tens\cdots \tens \gamma_n$. We proceed by induction over $n$: for $n=0$ the property
  is satisfied because $\alpha\subset 1_\Cat$ $\Leftrightarrow$ $\alpha = 1_\Cat$. Assume now
  that the property if satisfied for a given $n\geq 0$ and consider an $\alpha\in\Irr\Cat$ such
  that $|\alpha| = n+1$. By definition of $\Gg$ there exist $\beta\in \Vv$ and $\gamma\in\Dir$
  such that $|\beta| = n$ and $\alpha\subset \beta\tens\gamma$, and the induction hypothesis
  for $\beta$ gives the desired inclusion $\alpha\subset \gamma_1\tens \cdots \tens
  \gamma_n\tens \gamma$.  Assume conversely that $\alpha\subset \gamma_1\tens \cdots \tens
  \gamma_{n+1}$ and let $(\beta_k)$ be a maximal orthogonal family of irreducible subobjects of
  $\gamma_1\tens \cdots \tens \gamma_n$. Because $\alpha$ is irreducible the inclusion
  $\alpha\subset \dirsum (\beta_k\tens \gamma_{n+1})$ implies that $\alpha\subset \beta_k\tens
  \gamma_{n+1}$ for some $k$. By induction hypothesis one has $|\beta_k| \leq n$, hence
  $|\alpha| \leq n+1$.
  
  \ref{item:strict_repr}. Recall that the endpoints map $e$ induces a morphism from $\GG$
  onto $\Gg$, the latter one being a tree.  Therefore $\GG$ is a tree \iff $e$ is
  injective, and it is enough to check it on the ascending orientation $\Ee_+ \subset
  \Ee$: this leads to the condition that the subobjects $(\alpha\tens\gamma)_+$, for a
  given $\alpha$, should have pairwise different subobjects without multiplicity. The tree
  $\GG$ is then directional with respect to the origin $1_\Cat$ \iff the corepresentations
  $(\alpha\tens\gamma)_+$ have at most one subobject.
  
  \ref{item:always_asc}. and \ref{item:dim_grows}. We proceed again by induction on the
  distance to the origin: let $(\alpha,\beta)$ be an ascending edge with direction
  $\gamma$ and assume that $\dim\beta \geq \dim\alpha$, with equality \iff $\dim\gamma >
  1$. Take $\gamma'\in\Dir$, the assumption on $\GG$ shows that $(\beta\tens\gamma') =
  (\beta\tens\gamma')_+$ or $(\beta\tens\gamma') = \alpha \oplus (\beta\tens \gamma')_+$.
  In the first case, which can only happen when $(\beta\tens\gamma')_+ \neq 0$, one has
  clearly $\dim (\beta\tens\gamma')_+ \geq \dim\beta$ with equality \iff $\dim\gamma' =
  1$. On the other hand, we are in the second case \iff $\gamma' = \bar\gamma$, because of
  the equivalence $\alpha \subset \beta\tens\gamma' \Leftrightarrow \beta \subset
  \alpha\tens \bar\gamma'$. Moreover one has then $(\beta\tens\gamma')_+ = 0$ \iff
  $\dim\beta \dim\gamma' = \dim\alpha$, which is equivalent to $\dim\gamma' = 1$ by
  induction hypothesis. If on the contrary $\dim\gamma' = \dim\gamma > 1$, the strict case
  of the induction hypothesis gives
  \begin{displaymath}
    \dim (\beta\tens\gamma')_+ = \dim\beta \dim\gamma' - \dim\alpha \geq 2 \dim\beta -
    \dim\alpha > \dim\beta \text{.} 
  \end{displaymath}
\end{dem}

\begin{prp}\label{prp:free}
  Let $S$ be a full Woronowicz \Cst algebra and $p_1$ a central projection of $\hat S$ such
  that $Up_1U = p_1$ and $p_0 p_1 = 0$. If the full classical Cayley graph $\GG$ is a directional
  tree, then
  \begin{itemize}
  \item $S$ is a free product of a finite number of free Woronowicz \Cst algebras $A_o(Q_i)$
    and $A_u(Q'_j)$, with $Q_i\bar Q_i\in\CC\id$ and $Q'_j$ invertible,
  \item $p_1$ is the sum of the central supports of the respective fundamental
    corepresentations of these Woronowicz \Cst algebras.
  \end{itemize}
  Conversely the full classical Cayley graph $\GG$ of any such pair $(S,p_1)$ is a directional
  tree.
\end{prp}

\begin{dem}
  The classical graph $\Gg$ being a tree, the set $\Irr\Cat$ of its vertices lies in
  one-to-one correspondence with the set of paths without half-turns starting from the
  origin $1_\Cat$. Because $\Gg$ is isomorphic to the full graph $\GG$, these paths are
  characterized by the finite sequences of the directions they follow.  Finally,
  Lemma~\ref{lem:strict} shows that the finite sequences $(\gamma_i)$ of elements of
  $\Dir$ that arise in such a way are exactly the ones that fulfill the condition
  $\gamma_{i+1} \neq \bar \gamma_i$ or $\mathrm{dim}~ \gamma_{i+1}>1$ for each $i$.
  
  For every pair $\{\gamma, \bar\gamma\} \subset \Dir$ with $\gamma = \bar\gamma$ (resp.
  $\gamma\neq \bar\gamma$), the universal property of free quantum groups gives a Hopf
  homomorphism from some $A_o(Q)$ (resp.  $A_u(Q)$) onto $S$, where $Q$ is a matrix such
  that $Q\bar Q \in \CC\id$ (resp. is invertible). By universality of free products, one
  obtains then a surjective Hopf homomorphism $\Phi : F\to S$, where $F$ is some finite
  free product of free quantum groups. By definition, for each factor $A_o(Q)$, $A_u(Q)
  \subset F$ the fundamental corepresentation $U$ and its conjugate are mapped by
  $\id\tens\Phi$ onto the corresponding pair $\{\gamma, \bar\gamma\} \subset \Dir$.
  
  On the other hand, the starting remarks on the structure of $\Gg$ show that $\Irr\Cat$
  is the monoid generated by $\Dir$ and the relations $\{\gamma\bar\gamma = \bar\gamma
  \gamma = 1 ~|~ \gamma\in\Dir,~\dim\gamma = 1\}$. Hence $\Phi$ induces a bijection
  between $\Irr\Cat$ and the set $\Irr\mathcal{F}$ of irreducible corepresentations of $F$
  (up to equivalence) --- see \cite{Wang:freeprod, Banic:O(n)_cras, Banic:U(n)} for the
  description of $\Irr\mathcal{F}$ and notice that $A_o(Q)$ and $A_u(Q)$ are respectively
  isomorphic to $C^*(\ZZ/2\ZZ)$ and $C^*(\ZZ)$ when $\dim Q=1$. This proves, using
  \cite{Woro:dual} and the fact that we are dealing with full Woronowicz \Cst algebras,
  that $\Phi$ is injective.  The statement that $\GG$ is a directional tree for any free
  product of free quantum groups follows easily from the above mentioned description of
  $\Irr\mathcal{F}$.
\end{dem}

\begin{figure}[h!]    
  \begin{center}
    \begin{picture}(0,0)%
      \includegraphics{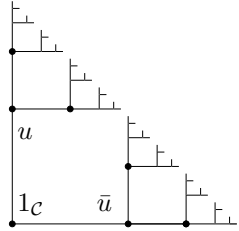}%
    \end{picture}%
    \setlength{\unitlength}{2368sp}%
    \begingroup\makeatletter\ifx\SetFigFont\undefined%
    \gdef\SetFigFont#1#2#3#4#5{%
      \reset@font\fontsize{#1}{#2pt}%
      \fontfamily{#3}\fontseries{#4}\fontshape{#5}%
      \selectfont}%
    \fi\endgroup%
    \begin{picture}(2375,2374)(5663,-4898)
      \put(6601,-4711){\makebox(0,0)[lb]{\smash{\SetFigFont{10}{13.2}{\familydefault}{\mddefault}{\updefault}{\color[rgb]{0,0,0}$\bar u$}%
          }}}
      \put(5776,-4711){\makebox(0,0)[lb]{\smash{\SetFigFont{10}{13.2}{\familydefault}{\mddefault}{\updefault}{\color[rgb]{0,0,0}$1_{\mathcal C}$}%
          }}}
      \put(5776,-3961){\makebox(0,0)[lb]{\smash{\SetFigFont{10}{13.2}{\familydefault}{\mddefault}{\updefault}{\color[rgb]{0,0,0}$u$}%
          }}}
    \end{picture}
  \end{center}
  \caption{classical Cayley graph of the unitary free quantum group}
  \label{fig:cayley}
\end{figure}

\begin{rque}{Example}(free quantum groups) Let us picture the simplest
  cases of Proposition~\ref{prp:free}. When $\dim Q>1$ and $Q\bar Q\in\CC\id$, the
  classical Cayley graph of $A_o(Q)$ endowed with its fundamental corepresentation is the
  half line with vertices at the integers. When $\dim Q>1$ and $Q$ is invertible, the
  classical Cayley graph of $A_u(Q)$ is drawn in Figure~\ref{fig:cayley}.
\end{rque}

\begin{prp} \label{prp:but_inj}
  Let $S$ be a Woronowicz \Cst algebra and $p_1$ a central projection of $\hat S$ such that
  $Up_1U = p_1$ and $p_0 p_1 = 0$. Assume that the classical Cayley graph $\GG$ of $(S, p_1)$
  is a directional tree.  Then the restriction of $\T$ to $K_{\ss++}$ is injective and its image
  is $(1 - p_0)H$.
\end{prp}

\begin{dem}
  Let $\alpha\in\Irr\Cat$ and $\gamma\in\Dir$ be such that $|\alpha| = n$ and
  $(\alpha\tens \gamma)_+ \neq 0$. The subspace $(p_\alpha \tens p_\gamma)K$ is
  irreducible with respect to the representation $\hat\pi_4$ of $\hat S^{\tensexp 4}$, and
  equivalent to $(\alpha\tens\gamma) \tens \overline{(\alpha\tens\gamma)}$.  By
  definition, $p_{\ss\jok+}(p_\alpha\tens p_\gamma) = \hat\delta(p_{n+1}) (p_\alpha\tens
  p_\gamma)$, so that $p_{\ss\jok+}(p_\alpha\tens p_\gamma)K$ is equivalent to
  $(\alpha\tens\gamma)_+ \tens \overline{(\alpha\tens\gamma)}$ with respect to the
  representation $\hat\pi_4\rond (\hat\delta\tens\id\tens\id)$ of $\hat S^{\tensexp 3}$.
  Similarly, and thanks to the first point of Lemma~\ref{lem:strict}, the subspace
  $p_{\ss+\jok}(p_\alpha\tens p_\gamma)K$ is equivalent to $(\alpha\tens\gamma) \tens
  \overline{(\alpha\tens\gamma)_+}$ with respect to the representation $\hat\pi_4\rond
  (\id\tens\id\tens\hat\delta)$. Finally, $p_{\ss++}(p_\alpha\tens p_\gamma)K$ is
  equivalent to $(\alpha\tens\gamma)_+ \tens \overline{(\alpha\tens\gamma)_+}$ for the
  representation $\hat\pi_4 \rond (\hat\delta\tens\hat\delta)$ of $\hat S\tens \hat S$,
  and therefore irreducible by hypothesis.
  
  Recall now from Proposition~\ref{prp:commutdual} that $\T$ intertwines $\hat\pi_4 \rond
  (\hat\delta\tens\hat\delta)$ and $\hat\pi_2$. Hence the restriction of $\T$ to
  $(p_\alpha\tens p_\gamma)K_{\ss++}$ is a multiple of an isometry, and one can compute
  the corresponding norm by considering the image of particular vectors, for instance
  characters. One has by Proposition~\ref{prp:endpoints}
  \begin{eqnarray*}
    \T p_{\ss++}(\chi_\alpha \tens \chi_\gamma) &=& 
    \T p_{\ss\jok+}(\chi_\alpha \tens \chi_\gamma) = 
    \T \hat\delta(p_{n+1})(\chi_\alpha \tens \chi_\gamma) \\
    &=& p_{n+1} \T (\chi_\alpha \tens \chi_\gamma) = 
    p_{n+1} (\chi_{\alpha\tensexp\gamma}) \\
    &=&  p_{n+1} (\chi_{(\alpha\tensexp\gamma)_-} 
    + \chi_{(\alpha\tensexp\gamma)_+}) = \chi_{(\alpha\tensexp\gamma)_+} \text{.}
  \end{eqnarray*}
  The norm in $H$ of the character of an irreducible corepresentation equals $1$ (cf
  \cite{Woro:cmp}, th.~5.8), so that one eventually gets the following lower bound for the
  norm of $\T p_{\ss++}(p_\alpha\tens p_\gamma)$ :
  \begin{equation}
    || \T p_{\ss++}(p_\alpha\tens p_\gamma) || = 
    \frac{||\chi_{(\alpha\tensexp\gamma)_+}||}
    {||p_{\ss++}(\chi_\alpha\tens\chi_\gamma)||}   
    \geq \frac {||\chi_{(\alpha\tensexp\gamma)_+}||} 
    {||\chi_\alpha\tens\chi_\gamma||} = 1 \text{.} 
    \label{eq:norme_but}    
  \end{equation}
  
  To conclude, let us remark that $\T p_{\ss++}$ maps the respective orthogonal subspaces
  $p_{\ss++}(p_\alpha\tens p_\gamma)K$ onto the subspaces $p_{(\alpha\tensexp\gamma)_+}H$, which
  are pairwise different by hypothesis, hence orthogonal, and whose sum equals $(1-p_0)H$. The
  operator $\T$ is therefore injective and has dense image in $(1-p_0)H$, but this image is
  closed by~(\ref{eq:norme_but}).
\end{dem}

\begin{rque}{Remark} \label{rque:ext_target}
  We will need in Section~\ref{section:AO} to have a slightly more general and more
  precise result than~(\ref{eq:norme_but}). Let $\Hh$ be the algebraic direct sum of the
  subspaces $p_kH$, and let $\Te$ be the operator defined on $\Hh\tens\Hh$ in the same way
  as $\T$, that is, by Proposition~\ref{prp:endpoints}, coming from the multiplication of
  $S$. Let us also extend $p_{\ss\jok+}$ to $\Hh\tens\Hh$ by putting $\Pe_{\ss\jok+} =
  \sum \hat\delta(p_{n+k}) (p_n\tens p_k)$: we have then $\Te \Pe_{\ss\jok+} (\id\tens
  p_1) = \T p_{\ss\jok+} = \T p_{\ss++}$. The first arguments of the preceding proof are
  still valid if one replaces $\gamma$ with any $\beta\in \Irr\Cat$: as a matter of fact
  the hypothesis implies that $\alpha \tens \beta$ contains at most one subobject $\delta$
  with $|\delta| = |\alpha| + |\beta|$. If such a $\delta$ exists one gets the following
  generalization of~(\ref{eq:norme_but}):
  \begin{displaymath}
    || \Te\Pe_{\ss\jok+}(p_\alpha\tens p_\beta) || = 
    ||\hat\delta(p_\delta) (\chi_\alpha\tens\chi_\beta)||^{-1} \text{.}
  \end{displaymath}
  Let us remark that $\chi_\alpha$ (resp. $t_\alpha(1)$) generates the invariant line of
  $p_\alpha H$ (resp. $H_\alpha\tens H_{\bar\alpha}$) with respect to the action of $\hat
  S$. Moreover one has $||\chi_\alpha|| = 1$ and $||t_\alpha(1)|| = \sqrt{M_\alpha}$, so
  that $\chi_\alpha$ in fact corresponds to $t_\alpha(1) / \sqrt{M_\alpha}$ in the
  isomorphism $p_\alpha H \simeq H_\alpha\tens H_{\bar\alpha}$, up to a phase factor.
  Consequently, in the isomorphism $p_\alpha H\tens p_\beta H \simeq$ $H_\alpha\tens
  H_\beta \tens H_{\bar\beta} \tens H_{\bar\alpha}$ the vector $\hat\delta(p_{\delta})
  (\chi_\alpha\tens\chi_\beta)$ corresponds to
  \begin{displaymath}
    (\hat\delta(p_{\delta})\tens\id\tens\id) (t_{\alpha\tensexp\beta}(1)) / 
    \sqrt{M_\alpha M_\beta} = t_{\delta}(1) / \sqrt{M_\alpha M_\beta} \text{,}
  \end{displaymath}
  if one isometrically identifies $H_\delta$ with the equivalent subspace of
  $H_\alpha\tens H_\beta$. We therefore get the following exact formula, from
  which~(\ref{eq:norme_but}) can be recovered by noticing that $M_\delta \leq
  M_\alpha\tens M_\beta$:
  \begin{displaymath}
    || \Te\Pe_{\ss\jok+}(p_\alpha\tens p_\beta) || = 
    \sqrt {\frac {M_\alpha M_\beta} {M_\delta}} \text{.} 
  \end{displaymath}
\end{rque}

\section{Geometric edges}
\label{section:geometric}

In this section we will study the Hilbert space $K_g = \Ker (\Theta + \id)$ when the
classical Cayley graph $\GG$ is a directional tree. We consider
Proposition~\ref{prp:but_inj} as an evidence that $K_{\ss++}$ provides a good notion of
``quantum ascending edges'', and we would similarly like to know whether $K_g$ provides a
good notion of ``quantum geometric edges''. By this we mean that there should be exactly
one geometric edge for each ascending edge, which can be more rigorously expressed in the
hilbertian framework by the fact that the restriction $p_{\ss++} : K_g \to K_{\ss++}$
should be invertible.

Of course the study of $K_g = \Ker (\Theta + \id)$ is closely related to the problem of
the non-involutivity of the reversing operator $\Theta$. The next proposition provides a
``weak involutivity'' property which we will use for the proof of Theorem~\ref{thm:proj},
as well as a technical corollary obtained in Lemma~\ref{lem:def_W}.  Notice that
$(p_{\ss++} + p_{\ss--})K$ behaves as a subspace of ``quasi-classical'' quantum edges in
this regard.

\begin{prp} \label{prp:sous_invol}
  Let $S$ be a Woronowicz \Cst algebra and $p_1$ a central projection of $\hat S$ such
  that $Up_1U = p_1$ and $p_0 p_1 = 0$. Assume that the classical Cayley graph $\GG$ of
  $(S, p_1)$ is a directional tree.  Then we have, for all $n\in\NN$:
  \begin{displaymath}
    (p_{\ss++} + p_{\ss--}) \Theta^n (p_{\ss++} + p_{\ss--}) = 
    (p_{\ss++} + p_{\ss--}) \Theta^{-n} (p_{\ss++} + p_{\ss--})
    \text{.} 
  \end{displaymath}
\end{prp}

\begin{dem}
  Inserting $\id = p_{\ss++} + p_{\ss+-} + p_{\ss-+} + p_{\ss--}$ between the occurrences
  of $\Theta^{\pm 1}$ in $\Theta^{\pm n}$ and developing, the statement of the theorem
  becomes an equality between two sums of terms looking like $p_{\epsilon'_0,\epsilon'_0}
  \Theta^{\pm 1} p_{\epsilon_1,\epsilon'_1} \Theta^{\pm 1} \cdots \Theta^{\pm 1}
  p_{\epsilon_n,\epsilon_n}$.  We will in fact prove that these terms are pairwise equal:
  for $\epsilon_i$, $\epsilon'_i \in \{+,-\}$ with $i \in \inter{0}{n}$, one has
  \begin{eqnarray} \nonumber
    && p_{\epsilon'_0,\epsilon'_0} \Theta p_{\epsilon_1,\epsilon'_1} 
    \Theta \cdots \Theta p_{\epsilon_{n-1},\epsilon'_{n-1}} \Theta
    p_{\epsilon_n,\epsilon_n} = \\ \label{eq:sous_invol_detail}
    && \makebox[2cm]{} = p_{\epsilon'_0,\epsilon'_0}
    \Theta^{-1} p_{\epsilon'_1,\epsilon_1} \Theta^{-1} \cdots
    \Theta^{-1} p_{\epsilon'_{n-1},\epsilon_{n-1}} \Theta^{-1}
    p_{\epsilon_n,\epsilon_n} \text{.} 
  \end{eqnarray}
  
  Let us proceed by induction over $n \in \NN$, calling ``rank $0$'' the trivial equality
  $p_{\epsilon',\epsilon'} p_{\epsilon,\epsilon}= p_{\epsilon',\epsilon'}
  p_{\epsilon,\epsilon}$. Choose $n\geq 1$. As a first step, assume that there exists $k
  \in \inter{1}{n-1}$ such that $\epsilon_k = \epsilon'_k$. Then the conclusion results
  straightforwardly from two applications of the induction hypothesis at ranks $k$ and
  $n-k$, with $(\epsilon_i,\epsilon'_i)_{0\leq i\leq k}$ and
  $(\epsilon_i,\epsilon'_i)_{k\leq i\leq n}$ respectively.
  
  We assume now that $\epsilon_i = - \epsilon'_i$ for each $i$. If one side
  of~(\ref{eq:sous_invol_detail}) is non-zero, we necessarily have $(\epsilon_i,
  \epsilon'_i) = (-\epsilon'_0,-\epsilon_n)$ for all indices $i$: as a matter of fact,
  Proposition~\ref{prp:orient} shows that the equalities $\epsilon_{i+1}= -\epsilon'_i$
  are required for the products in~(\ref{eq:sous_invol_detail}) not to vanish. In
  particular, we have then $\epsilon'_0 = -\epsilon_n$. This proves that the equalities
  from the first step are sufficient to get the identities $p_{\ss--}\Theta^n p_{\ss--} =
  p_{\ss--}\Theta^{-n} p_{\ss--}$ and $p_{\ss++}\Theta^n p_{\ss++} = p_{\ss++}\Theta^{-n}
  p_{\ss++}$.  Moreover, taking the adjoint allows to switch from $\epsilon'_0 = -1$ to
  $\epsilon'_0 = 1$, so that it only remains to prove the equality
  \begin{displaymath}
    p_{\ss++}\Theta p_{\ss-+}\Theta \cdots
    \Theta p_{\ss-+}\Theta p_{\ss--} = 
    p_{\ss++} \Theta^{-1} p_{\ss+-}\Theta^{-1}  
    \cdots\Theta^{-1}  p_{\ss+-}\Theta^{-1} p_{\ss--} \text{.} 
  \end{displaymath}
  
  By adding terms from the first step we rather focus on the following equivalent
  equality:
  \begin{displaymath}
    p_{\ss++} \Theta^n p_{\ss--} + p_{\ss--}\Theta^n p_{\ss--}
    \stackrel{\text{?}}{\ds =} 
    p_{\ss++} \Theta^{-n} p_{\ss--} + p_{\ss--}\Theta^{-n} p_{\ss--}
    \text{.}
  \end{displaymath}
  Using the fact from Proposition~\ref{prp:but_inj} that the target operator $\T$ is
  injective on $K_{\ss++}$, we can compose on the left by $\T$ and use
  Proposition~\ref{prp:orient} to get another equivalent equality: $\T \Theta^n p_{\ss--}
  = \T \Theta^{-n} p_{\ss--}$. But this is true since we have, from
  Proposition~\ref{prp:endpoints} and the definition of $\T$: $\T\Theta^2 = \S\Theta = \T$,
  hence $\T \Theta^{2k} = \T \Theta^{-2k} = \T$ and $\T\Theta = \T\Theta^{-1}$.
\end{dem}

\begin{lem} \label{lem:def_W}
  Let $S$ be a Woronowicz \Cst algebra and $p_1$ a central projection of $\hat S$ such
  that $Up_1U = p_1$ and $p_0 p_1 = 0$. Assume that the classical Cayley graph $\GG$ of
  $(S, p_1)$ is a directional tree. Then there exists a unique unitary operator $W :
  K_{\ss+-} \to K_{\ss-+}$ such that
  \begin{displaymath}
    \forall~ k\in\NN~~
    W (p_{\ss+-}\Theta)^k p_{\ss++} = 
    (p_{\ss-+}\Theta^{-1})^k p_{\ss++} \text{.}
  \end{displaymath}
  Moreover we have $W p_{\ss+-}\Theta = p_{\ss-+}\Theta^{-1} W$ and $p_{\ss--}\Theta =
  p_{\ss--}\Theta^{-1} W$ on $K_{\ss+-}$.
\end{lem}

\begin{dem}
  Let $X$ (respectively $X'$) be the operator from $K_{\ss++}\tens \ell_2(\NN)$ to $K_{\ss+-}$
  (resp.  $K_{\ss-+}$) defined by $X(\xi\tens e_k) = 2^{-k} (p_{\ss+-}\Theta)^k \xi$
  (resp. $X'(\xi\tens e_k) =$ $2^{-k} (p_{\ss-+}\Theta^{-1})^k \xi$). Thanks to the
  coefficients $2^{-k}$, the operators $X$ and $X'$ are bounded, and it is easy to see
  that their adjoints are resp. given by
  \begin{displaymath}
    X^* = \sum 2^{-k} T_k p_{\ss++} (\Theta^{-1}p_{\ss+-})^k 
    \text{~~~and~~~} X^{\prime *} = 
    \sum 2^{-k} T_k p_{\ss++} (\Theta p_{\ss-+})^k \text{,}
  \end{displaymath}
  where we put $T_k(\xi) = \xi\tens e_k$ for any $\xi\in K_{\ss++}$.  Let $\zeta$ be an
  element of $\Ker X^*$, for every $k$ and $n$ we have $p_{\ss++} (\Theta ^{-1}
  p_{\ss+-})^k (p_n\tens\id) \zeta = 0$.  In particular $(\Theta^{-1} p_{\ss+-})^n
  (p_n\tens p_1) \zeta$ vanishes: by Proposition~\ref{prp:orient} it is an element of
  $(p_0\tens\id) K$, which is contained in $K_{\ss++}$. By a finite descending induction
  on $k\in \inter{0}{n}$, we deduce that
  \begin{displaymath}
    (\Theta^{-1} p_{\ss+-})^k (p_n\tens p_1)\zeta = 
    p_{\ss++} (\Theta^{-1} p_{\ss+-})^k (p_n\tens p_1)\zeta + 
    \Theta (\Theta^{-1} p_{\ss+-})^{k+1}(p_n\tens p_1)\zeta 
  \end{displaymath}
  vanishes, and in particular $(p_n\tens\id)\zeta = 0$ for any $n$.  Hence $X^*$ is
  injective and $X$ has dense image. In the same way, $X'$ has dense image.
  
  To prove the existence and the uniqueness of $W$, which is characterized by the identity $WX
  = X'$, it is therefore enough to show that $||X\eta|| = ||X'\eta||$ for any $\eta\in
  K_{\ss++}\tens \ell_2(\NN)$, or as well, that $X^*X = X^{\prime *}X'$. We will work on each
  subspace $K_{\ss++}\tens e_i$ separately, and we are thus led to prove for every $k$ and $l$
  the equality
  \begin{equation} \label{eq:dem_def_W}
    p_{\ss++}(\Theta^{-1}p_{\ss+-})^l(p_{\ss+-}\Theta)^kp_{\ss++} =
    p_{\ss++}(\Theta p_{\ss-+})^l(p_{\ss-+}\Theta^{-1})^kp_{\ss++} 
    \text{,}
  \end{equation}
  which can also be written
  \begin{eqnarray*}
    && p_{\ss++}\Theta^{-1} p_{\ss+-}\cdots p_{\ss+-}\Theta^{-1}
    (1 - p_{\ss--}) \Theta p_{\ss+-}\cdots p_{\ss+-} \Theta p_{\ss++}
    = \\ && \makebox[2cm]{} = 
    p_{\ss++}\Theta p_{\ss-+}\cdots p_{\ss-+}\Theta (1 - p_{\ss--}) 
    \Theta^{-1} p_{\ss-+}\cdots p_{\ss-+}\Theta^{-1} p_{\ss++}
    \text{.} 
  \end{eqnarray*}
  We proceed by induction on $\min(k,l)$ and distribute $(1-p_{\ss--})$ on both sides: the
  terms coming from $p_{\ss--}$ are equal thanks to Equation~(\ref{eq:sous_invol_detail})
  of Proposition~\ref{prp:sous_invol}, and the terms coming from $1$ are equal by
  induction hypothesis. When $kl=0$ but $(k,l)\neq (0,0)$, both sides
  of~(\ref{eq:dem_def_W}) vanish, and when $k=l=0$, (\ref{eq:dem_def_W}) is trivial.
  
  Because $X$ has dense image, it suffices to check the equalities $W p_{\ss+-}\Theta =
  p_{\ss-+}\Theta^{-1} W$ and $p_{\ss--}\Theta = p_{\ss--}\Theta^{-1} W$ on the image of
  $(p_{\ss+-}\Theta)^k p_{\ss++}$, for every $k$. The first one follows immediately from
  the definition of $W$:
  \begin{eqnarray*}
    (Wp_{\ss+-}\Theta)~ (p_{\ss+-}\Theta)^k p_{\ss++} &=& 
    W (p_{\ss+-}\Theta)^{k+1} p_{\ss++} 
    = (p_{\ss-+}\Theta^{-1})^{k+1} p_{\ss++} \text{~~~and} \\
    (p_{\ss-+}\Theta^{-1} W)~ (p_{\ss+-}\Theta)^k p_{\ss++} &=&
    p_{\ss-+}\Theta^{-1} (p_{\ss-+}\Theta^{-1})^k p_{\ss++} 
    = (p_{\ss-+}\Theta^{-1})^{k+1} p_{\ss++} \text{.}
  \end{eqnarray*}
  For the second one, we use furthermore Equation~(\ref{eq:sous_invol_detail}) from
  Proposition~\ref{prp:sous_invol}:
  \begin{eqnarray*}
    (p_{\ss--}\Theta^{-1} W)~ (p_{\ss+-}\Theta)^k p_{\ss++} &=&
    p_{\ss--}\Theta^{-1} (p_{\ss-+}\Theta^{-1})^k p_{\ss++} \\ &=& 
    (p_{\ss--}\Theta)~ (p_{\ss+-}\Theta)^k p_{\ss++} \text{.}
  \end{eqnarray*} 
\end{dem}

\begin{thm} \label{thm:proj}
  Let $S$ be a Woronowicz \Cst algebra and $p_1$ a central projection of $\hat S$ such
  that $Up_1U = p_1$ and $p_0 p_1 = 0$. Assume that the classical Cayley graph $\GG$ of
  $(S, p_1)$ is a directional tree. Then the orthogonal projection from $K_g$ to $K_{\ss++}$
  is injective and its image is given by
  \begin{displaymath}
    p_{\ss++} K_g = \{ \zeta \in K_{\ss++} ~|~ 
    \exists~\eta \in K_{\ss+-}~~ 
    (\id + p_{\ss+-}\Theta)(\eta) = p_{\ss+-}\Theta\zeta \} \text{.} 
  \end{displaymath}
\end{thm}

\begin{dem}
  By definition, a vector $\xi\in K$ lies in $K_g$ \iff $\Theta(\xi) = -\xi$, which we
  split in two equations: $p_{\ss\jok-}\xi = - \Theta p_{\ss+\jok}\xi$ and
  $p_{\ss-\jok}\xi = - \Theta^{-1} p_{\ss\jok+}\xi$. Let us first analyze these conditions
  with respect to the decomposition $K = \dirsum (p_k\tens\id) K$, using
  Proposition~\ref{prp:orient}:
  \begin{eqnarray*}
    \forall~ n\in\NN ~~ p_{\ss\jok-}(p_n\tens\id)\xi &=& 
    - \Theta p_{\ss+\jok} (p_{n-1}\tens\id)\xi \text{~~~and} \\
    \forall~ n\in\NN ~~ p_{\ss-\jok}(p_n\tens\id)\xi &=&
    - \Theta^{-1}p_{\ss\jok+} (p_{n-1}\tens\id)\xi \text{.}
  \end{eqnarray*}
  If $p_{\ss++}\xi = 0$, this gives a linear induction equation for $((p_n\tens\id)
  \xi)_n$, and since $(p_0\tens\id)\xi = p_{\ss++}(p_0\tens\id)\xi = 0$ the whole sequence
  vanishes. Hence the restriction of $p_{\ss++}$ to $K_g$ is injective.
  
  Now we use the decomposition $\id = p_{\ss++} + p_{\ss+-} + p_{\ss-+} + p_{\ss--}$ to
  get a new system equivalent to the conditions $p_{\ss\jok-} \xi = - \Theta p_{\ss+\jok}
  \xi$ and $p_{\ss-\jok} \xi = - \Theta^{-1} p_{\ss\jok+}$, which characterize vectors in
  $K_g$:
  \newcounter{dem_proj_syst}\setcounter{dem_proj_syst}{\value{equation}}%
  \begin{eqnarray} 
    p_{\ss+-}\xi &=& - p_{\ss+-}\Theta p_{\ss++}\xi -
    p_{\ss+-}\Theta p_{\ss+-}\xi \label{eq:dem_proj_syst_cond} \\
    p_{\ss--}\xi &=& - p_{\ss--}\Theta p_{\ss++}\xi -
    p_{\ss--}\Theta p_{\ss+-}\xi \\
    p_{\ss--}\xi &=& - p_{\ss--}\Theta^{-1} p_{\ss++}\xi -
    p_{\ss--}\Theta^{-1}p_{\ss-+}\xi \\
    p_{\ss-+}\xi &=& - p_{\ss-+}\Theta^{-1} p_{\ss++}\xi -
    p_{\ss-+}\Theta^{-1}p_{\ss-+}\xi \text{.} 
  \end{eqnarray}
  Let $\zeta\in K_{\ss++}$ be as in the statement of the theorem: there exists $\eta\in
  K_{\ss+-}$ such that $(\id + p_{\ss+-}\Theta)(\eta) = p_{\ss+-}\Theta\zeta$. Put $\xi =
  \zeta - \eta - W\eta + p_{\ss--}\Theta(\eta - \zeta)$. In this case, the above system
  can be written in the following way:
  \let\theorigeq\theequation \renewcommand\theequation{\theorigeq'}%
  \newcounter{saveeq}\setcounter{saveeq}{\value{equation}}%
  \setcounter{equation}{\value{dem_proj_syst}}%
  \begin{eqnarray}
    - \eta &=& - p_{\ss+-}\Theta \zeta + p_{\ss+-}\Theta \eta 
    \label{eq:dem_proj_syst_+-} \\ 
    p_{\ss--}\Theta(\eta - \zeta) &=& - p_{\ss--}\Theta\zeta 
    + p_{\ss--}\Theta \eta \label{eq:dem_proj_syst_--dir} \\
    p_{\ss--}\Theta(\eta - \zeta) &=& - p_{\ss--}\Theta^{-1} \zeta  
    + p_{\ss--}\Theta^{-1} W\eta \label{eq:dem_proj_syst_--inv} \\ 
    - W \eta &=& - p_{\ss-+}\Theta^{-1}\zeta 
    + p_{\ss-+}\Theta^{-1}W\eta \text{.} \label{eq:dem_proj_syst_-+} 
  \end{eqnarray}
  \let\theequation\theorigeq%
  \setcounter{equation}{\value{saveeq}}%
  We can notice that~(\ref{eq:dem_proj_syst_+-}) amounts to the hypothesis on $\zeta$ and
  $\eta$, whereas~(\ref{eq:dem_proj_syst_--dir}) is trivial.
  Proposition~\ref{prp:sous_invol} and Lemma~\ref{lem:def_W} show
  that~(\ref{eq:dem_proj_syst_--inv}) is always satisfied. Finally the hypothesis on
  $\zeta$ and $\eta$ yields $W\eta = W p_{\ss+-}\Theta\zeta - Wp_{\ss+-}\Theta\eta$,
  and~(\ref{eq:dem_proj_syst_-+}) follows then from Lemma~\ref{lem:def_W}.  Hence $\xi$
  lies in $K_g$ and $\zeta = p_{\ss++}\xi$ is in $p_{\ss++}K_g$. The reverse inclusion can
  easily be obtained from~(\ref{eq:dem_proj_syst_cond}): if $\zeta$ equals $p_{\ss++}\xi$
  with $\xi\in K_g$, we put $\eta = p_{\ss+-}\xi$ and the above mentioned equation reads
  then $(\id + p_{\ss+-} \Theta) (\eta) = p_{\ss+-} \Theta\zeta$, as already noticed.
\end{dem}

\section{Edges at infinity: the set}
\label{section:infinity_set}

The expression for $p_{\ss ++}K_g$ obtained in Theorem~\ref{thm:proj} is trivial in the
classical case because the projection $p_{\ss +-}$ vanishes then, but it has to be
analyzed in greater detail in the quantum case. More precisely, we need to understand the
interaction between $p_{\ss +-}$ and $\Theta$, and we will see that it can be described by
a purely quantum object: the space of ``edges at infinity'' $K_\infty$, that we introduce
in Definition~\ref{df:inf_edges}.

This definition bases on the simple remark that the operator $p_{\ss +-}\Theta p_{\ss +-}$
maps $(p_k\tens\id) K_{\ss+-}$ to $(p_{k+1}\tens\id) K_{\ss+-}$ by
Proposition~\ref{prp:orient}, and acts therefore as a right shift in the decomposition of
$K_{\ss+-}$ given by the distance to the origin in the classical Cayley graph.  It is then
very natural to introduce the associated inductive limit $K_{\infty}$.
Proposition~\ref{prp:RS} serves as a more precise motivation for this definition and shows
that the existence of $K_\infty$ is an obstruction to the surjectivity of $p_{\ss++} : K_g
\to K_{\ss++}$. Notice that in the classical case, the subspaces $(p_k\tens\id) K_{\ss+-}$
vanish, so that $K_\infty$ equals zero.

\begin{df} \label{df:inf_edges}
  Let $S$ be a Woronowicz \Cst algebra and $p_1$ a central projection of $\hat S$ such
  that $Up_1U = p_1$ and $p_0 p_1 = 0$. Assume that the classical Cayley graph associated
  with $(S,p_1)$ is a tree.
  \begin{enumerate}
  \item Put $r = - p_{\ss +-} \Theta p_{\ss+-}$, $s = p_{\ss+-} \Theta p_{\ss++}$ and
    define the inductive limit Hilbert space $K_{\infty} =
    {\ds\lim_{\longrightarrow}} ((p_k\tens\id) K_{\ss+-}, r)$.
  \item Let $R_k$ be the natural morphism from $(p_k\tens\id) K_{\ss+-}$ to $K_{\infty}$,
    and denote by $R$ the linear map $\sum_{k\geq 0} R_k$ defined on
    $\dirsum_{\mathrm{alg}} (p_k\tens\id) K_{\ss+-}$.
  \end{enumerate}
\end{df}

\begin{prp} \label{prp:RS}
  Let $S$ be a Woronowicz \Cst algebra and $p_1$ a central projection of $\hat S$ such
  that $Up_1U = p_1$ and $p_0 p_1 = 0$. Assume that the classical Cayley graph $\GG$ of
  $(S, p_1)$ is a directional tree.
  \begin{enumerate}
  \item \label{item:co_isom} The map $Rs$ extends to a co-isometry from $K_{\ss++}$ to
    $K_{\infty}$.
  \item \label{item:ker_geom} The subspace $p_{\ss++} K_g$ is contained in $\Ker Rs$.
    Moreover if the $R_k$ are injective one has, denoting by $p_{\geq k}$ the sum
    $\sum_{i\geq k} p_i\tens\id$:
    \begin{displaymath}
      p_{\ss++} K_g = \{ \zeta\in \Ker Rs ~|~ 
      (||R_k^{-1}Rs p_{\geq k}\zeta||)_k \in \ell^2(\NN) \} \text{.}
    \end{displaymath}
  \end{enumerate}
\end{prp}

\begin{dem}
  \ref{item:co_isom}. We start with a simple computation, using
  Proposition~\ref{prp:orient}:
  \begin{eqnarray*}
    rr^* + ss^* &=& p_{\ss+-}\Theta p_{\ss+-}\Theta^* p_{\ss+-} +
    p_{\ss+-} \Theta p_{\ss++} \Theta^* p_{\ss+-} \\
    &=& p_{\ss+-} \Theta p_{\ss+\jok} \Theta^* p_{\ss+-}
    = p_{\ss+-} \Theta \Theta^* p_{\ss+-} = \id_{K_{\ss+-}} \text{.}
  \end{eqnarray*}
  Notice that $R_0 = 0$ because $p_{\ss+-}(p_0\tens\id) = 0$, and that $R_{k+1}r = R_k$
  for any $k\in \NN$, by definition. We have then, denoting by $p_{\leq k}$ the sum
  $\sum_{i\leq k} p_i\tens\id$:
  \begin{eqnarray*}
    (Rs p_{\leq k}) (Rs p_{\leq k})^* &=& 
    \sum_{i=0}^{k-1} R_{i+1}s s^* R_{i+1}^* = 
    \sum_{i=0}^{k-1} R_{i+1} (1 - rr^*) R_{i+1}^* \\
    &=& \sum_{i=0}^{k-1} (R_{i+1} R_{i+1}^* - R_i R_i^*) = 
    R_k R_k^* \text{.}
  \end{eqnarray*}
  The maps $R_k$ being contractive, it follows that $Rs p_{\leq k}$ and $Rs$ itself are
  contractions. Because $p_{\leq k}$ converges to the identity in the $*$-strong topology,
  $(Rsp_{\leq k}) (Rsp_{\leq k})^*$ converges strongly to $(Rs)(Rs)^*$, and so it
  remains to show that $R_k R_k^*$ converges to the identity of $K_{\infty}$.  This is
  actually a general fact for contractive inductive limits: for any $l\geq k\geq 0$ and
  any $y\in (p_k\tens\id) K_{\ss+-}$, we have
  \begin{eqnarray*}
    ||R_lR_l^* (R_ky) - R_ky||^2 &\leq&
    ||R_l^*R_k y - r^{l-k} y||^2 \text{~~~and} \\ 
    ||R_l^*R_k y - r^{l-k} y||^2 &=& ||R^*_lR_k y||^2 -2\Re (R^*_lR_k y | 
    r^{l-k}y) + ||r^{l-k}y||^2 \\ &=& ||R^*_lR_k y||^2 - 2 ||R_k y||^2 +
    ||r^{l-k}y||^2 \\ &\leq& ||r^{l-k}y||^2 - ||R_k y||^2 \text{.}
  \end{eqnarray*}
  This upper bound tends to zero as $l$ goes to infinity, by definition of the norm of
  $K_\infty$. The union $\cup \Im R_k$ being dense in $K_{\infty}$, this proves that $R_lR_l^*
  \rightarrow_s \id$.
  
  \ref{item:ker_geom}. Let $\zeta \in p_{\ss++}K_g$: by Theorem~\ref{thm:proj}, there
  exists $\eta\in K_{\ss+-}$ such that $(1-r) \eta = s\zeta$. This can also be written
  \begin{eqnarray} \nonumber
    && \forall~k \in \NN^*~~ (p_k\tens\id) \eta = s (p_{k-1}\tens\id)
    \zeta + r (p_{k-1}\tens\id) \eta \\ \nonumber
    &\Longleftrightarrow& \forall~k \in \NN^*~~ (p_k\tens\id) \eta =
    \sum_{i=0}^{k-1} r^{k-i-1}s (p_{i}\tens\id) \zeta \\
    \label{eq:partial_sums} &\Longrightarrow& \forall~k \in \NN^*~~ 
    R_k(p_k\tens\id) \eta = Rs p_{\leq k-1} \zeta \text{.}  
  \end{eqnarray}
  The right-hand side of this equality converges to $Rs\zeta$ when $k$ goes to infinity,
  whereas the left-hand side tends to zero. Hence $p_{\ss++}K_g \subset \Ker Rs$.  Now, if
  the $R_k$ are injective, the implication leading to~(\ref{eq:partial_sums}) is an
  equivalence, so that a vector $\zeta\in K_{\ss++}$ is in $p_{\ss++} K_g$ \iff
  (\ref{eq:partial_sums}) defines a vector $\eta\in K_{\ss+-}$ \iff the orthogonal
  sequence $(R_k^{-1}Rs p_{\leq k-1}\zeta)_k$ is summable in $K_{\ss+-}$. Finally, we have 
  clearly $Rs p_{\leq k-1} \zeta = -Rs p_{\geq k} \zeta$ when $\zeta$ lies in $\Ker Rs$.
\end{dem}  

\bigskip

The rest of this section will be devoted to a more detailed study of $K_\infty$. We first
want to compute exactly the weights of the ``shift'' $p_{\ss+-}\Theta p_{\ss+-}$: this is
accomplished in Lemma~\ref{lem:lecture_theta} and relies on the technical results of
Section~\ref{section:corepr}. It is then easy to show that the maps $R_k$ are indeed
injective, and therefore that $K_\infty$ is infinite-dimensional in the quantum case.
Using the explicit result of Lemma~\ref{lem:lecture_theta} we are also able in
Theorem~\ref{thm:proj2} to make more precise the second statement of
Proposition~\ref{prp:RS}: it appears that $K_{\infty}$ is the only obstruction the
non-surjectivity of $p_{\ss++} : K_g \to K_{\ss++}$, except when the free product $(S,
\delta)$ under consideration contains one of the ``exceptional cases''
$A_o\big({0\atop-1}{1\atop0}\big)$, $A_o\big({0\atop 1}{1\atop0}\big)$ and
$A_u\big({1\atop0}{0\atop1}\big)$.

\bigskip

Let $(\gamma_1, \ldots, \gamma_k)$ be a finite sequence of directions $\gamma_i\in\Dir$.
There is at most one vertex $\alpha\in \Irr\Cat$ such that the geodesic from $1_\Cat$ to
$\alpha$ follows successively these directions: we will then put $\alpha = \gamma_1 \cdots
\gamma_k$. Now choose $\gamma\in\Dir$ and put $\gamma_{2l} = \bar\gamma$, $\gamma_{2l+1} =
\gamma$. We denote by $\alpha_k$ the vertex $\gamma \bar\gamma \cdots \gamma_k$, when it
exists. Lemma~\ref{lem:strict} shows that the set of values of $k$ is $\{0, 1\}$ when
$\dim \gamma=1$ and $\NN$ otherwise. We define in both cases the associated projection
$P_{\gamma} = \sum p_{\alpha_k} \tens p_{\bar\gamma_k}$. It is a central element of
$M(\hat S\tens \hat S)$, hence it commutes to the projections $p_{\ss\jok+}$ and
$p_{\ss+\jok}$. Moreover one has by Proposition~\ref{prp:commutdual}:
\begin{eqnarray*}
  \Theta P_{\gamma} &=& \sum \Theta (p_{\alpha_k} \tens p_{\bar\gamma_k}) 
  = \sum \hat\delta (p_{\alpha_k}) (1\tens p_{\gamma_k}) \Theta \\
  &=& \sum (p_{\ss\jok+} (p_{\alpha_{k-1}}\tens p_{\gamma_k}) +
  p_{\ss\jok-} (p_{\alpha_{k+1}}\tens p_{\gamma_k})) \Theta
  = P_\gamma \Theta \text{.}
\end{eqnarray*}

On the other hand, the projections $p_{\ss+\jok}$ and $p_{\ss\jok+}$ commute respectively
to the representations $\hat\pi_4\rond (\hat\delta\tens\id\tens\id)$ and $\hat\pi_4\rond
(\id\tens\id\tens\hat\delta)$ of $\hat S\tens\hat S\tens\hat S$ on $K$, by definition. In
particular they both commute to the representation $\hat\pi_4 \rond \hat\delta^3$ of $\hat
S$, as well as $\Theta$: see Proposition~\ref{prp:commutdual}. Hence $p_{\ss+\jok}$,
$p_{\ss\jok+}$ and $\Theta$ all commute to the projections $q_l =
\hat\pi_4\hat\delta^3(p_{2l})$. Let us recall in the case $\dim \gamma > 1$ that
$(p_{\alpha_k} \tens p_{\bar\gamma_k}) K_{\ss+-}$ is equivalent to $\alpha_{k-1} \tens
\bar\alpha_{k+1}$ with respect to $\hat\pi_4 \rond \hat\delta^3$, so that $q_l
(p_k\tens\id) P_\gamma K_{\ss+-}$ is non-zero \iff $l\in \inter{1}{k}$, and is then
irreducible and equivalent to $\alpha_{2l}$.

\begin{lem} \label{lem:lecture_theta}
  Let $\Theta$ be the reversing operator of a quantum Cayley tree, and choose $\gamma \in \Dir$
  with $\dim \gamma > 1$. We put $m_k = M_{\alpha_k}$ and $m_{-1} = 0$. Let $\epsilon_1$,
  $\epsilon'_1$, $\epsilon_2$, $\epsilon'_2 \in\{+, -\}$, with $\epsilon'_2 = -\epsilon_1$. For
  any $k \geq 1$ and $l \in \inter{1}{k}$ the operator $p_{\epsilon_2,\epsilon'_2} \Theta
  p_{\epsilon_1,\epsilon'_1}$ is a multiple of an isometry on $(p_{k-\epsilon_1}\tens\id)
  p_{\epsilon_1,\epsilon'_1} q_l P_\gamma K$ and
  \begin{displaymath}
    ||(p_k\tens p_1) p_{\epsilon_2,\epsilon'_2} \Theta 
    p_{\epsilon_1,\epsilon'_1}  q_l P_\gamma||= 
    \left\{ \begin{array}{ll}
        \sqrt{\ts\frac{m_lm_{l-1}}{m_{k}m_{k-1}}\ds} &\text{if  
          $\epsilon_1\epsilon'_1 \neq \epsilon_2\epsilon'_2$,} \\
        \sqrt{1-\ts\frac{m_lm_{l-1}}{m_{k}m_{k-1}}\ds} 
        &\text{if  $\epsilon_1\epsilon'_1 = \epsilon_2\epsilon'_2$.}
      \end{array} \right.
  \end{displaymath}
\end{lem} 

\begin{dem} 
  We can assume here that $S = A_u(Q)$ or $A_o(Q)$ because $\dim \gamma > 1$: see the
  proof of Proposition~\ref{prp:free}. We start with $p_{\ss+-}\Theta p_{\ss+-}$, by
  reorganizing the terms of the product and composing on the left by $\Theta^*$:
  \begin{displaymath}
    ||(p_k\tens\id) p_{\ss+-}\Theta p_{\ss+-} q_l P_\gamma|| =
    ||(\Theta^*p_{\ss+\jok}\Theta) p_{\ss\jok -} q_l ~ p_{\ss+\jok}
    (p_{\alpha_{k-1}} \tens p_{\gamma_k})|| \text{.}
  \end{displaymath}
  We know from the proof of Proposition~\ref{prp:but_inj} that the space
  $p_{\ss+\jok}(p_{\alpha_{k-1}}\tens p_{\gamma_k}) K$ is irreducible for the
  representation $\hat\pi_4 \rond (\id\tens\id\tens \hat\delta)$ of $\hat S\tens\hat
  S\tens\hat S$ and identifies with $\alpha_{k-1}\tens \gamma_k \tens \bar\alpha_{k}$. Let
  us study how $\Theta^*p_{\ss+\jok}\Theta$, $p_{\ss\jok -}$ and $q_l$ act in this
  identification.
  \begin{itemize}
  \item We have $p_{\ss+\jok} = \hat\pi_4(\hat\delta\tens\id\tens\id) (1\tens p)$, where
    $p = \sum \hat\delta(p_{n+1}) (p_1\tens p_{n})$. Lemma~\ref{prp:commutdual} shows 
    that $\Theta^*p_{\ss+\jok}\Theta = \hat\pi_4 (\id\tens\id\tens\hat\delta) (1\tens p)$,
    which hence acts on $\alpha_{k-1}\tens \gamma_k \tens \bar\alpha_{k}$ as $1\tens
    p$, ie as the projection onto $\alpha_{k-1}\tens \bar\alpha_{k+1}$.
  \item We know again from the proof of Proposition~\ref{prp:but_inj} that $p_{\ss\jok-}$
    acts in the identification like the projection of
    $\alpha_{k-1}\tens \gamma_k\tens \bar\alpha_{k}$ onto $\alpha_{k-2}\tens \bar\alpha_{k}$.
  \item Finally $q_l = \hat\pi_4 \hat\delta^3 (p_{2l})$ corresponds to the projection of
    $\alpha_{k-1}\tens \gamma_k \tens \alpha_k$ onto the sum of its subspaces that are
    equivalent to $\alpha_{2l}$.
  \end{itemize}
  Therefore Lemma~\ref{lem:angles_1} gives exactly the desired result for $||(p_k\tens\id)
  p_{\ss+-}\Theta p_{\ss+-}$ $q_l P_\gamma||$. We get then the norm of $(p_k\tens\id)
  p_{\ss+-}\Theta p_{\ss++} q_l P_\gamma$ by noticing that the sum of the squares of both
  norms equals $1$, and we proceed in the same way for the other cases.
\end{dem}

\begin{rque}{Remarks} \label{rque:subtrees}
  \begin{enumerate}
  \item \label{item:class_subtree} When $l=0$ --- and this includes the cases when $k=0$
    or $\dim \gamma = 1$---, we automatically have $p_{\ss+-} q_0 = p_{\ss-+} q_0 = 0$. In
    particular $p_{\ss++}\Theta p_{\ss--} q_0 = \Theta p_{\ss--} q_0$, and therefore
    Lemma~\ref{lem:lecture_theta} is replaced in this case by the statement that
    $p_{\ss++}\Theta p_{\ss--}$ and $p_{\ss--}\Theta p_{\ss++}$ are isometric on
    $p_{\ss--} q_0 K$ and $p_{\ss++} q_0 K$. In fact the subspace $q_0 K$, the analogous
    subspace $\hat\pi_2 \hat\delta(p_0) H$ and the corresponding restrictions of $\Theta$
    and $E$ are exactly the hilbertian objects associated to the classical Cayley graph
    $\Gg$.
  \item \label{item:gen_subtree} Lemma~\ref{lem:lecture_theta} only concerns the
    ``subtrees'' $P_\gamma K$, but this is enough to get results about the whole of $K$,
    thanks to a ``cut-and-paste'' process that we explain now. Let $\mathcal{I}$ be the
    set of ordered pairs $(\beta, \gamma) \in \Irr\Cat \times \Dir$ such that the last
    direction followed by the geodesic from $1_\Cat$ to $\beta$ is different from
    $\bar\gamma$ --- including $(1_\Cat, \gamma)$ for all $\gamma\in\Dir$. For such a
    $(\beta,\gamma)$ we denote by $\beta_k$ the vertices on the ascending path starting
    from $\beta$ and taking the directions $\gamma$, $\bar\gamma$, $\ldots\,$, and we call
    $P_{\beta,\gamma}$ the sum of the $p_{\beta_k} \tens p_{\bar\gamma_k}$.  Because the
    edges of the classical Cayley graph $\Gg$ are walked through once by exactly one of
    these paths, we see that $K$ is the orthogonal direct sum over $\mathcal{I}$ of the
    $P_{\beta,\gamma}K$. Notice that $P_\gamma = P_{1_\Cat, \gamma}$.
  
    Now we use the ``extended target operator'' $\Te : \Hh \tens H \to H$, ie the operator
    induced in the GNS construction of the Haar state by the multiplication of $S$. Take
    $(\beta, \gamma) \in \mathcal{I}$ and denote by $\alpha_k$ the objects constructed
    from $(1_\Cat, \gamma)$ as above. By definition of $\mathcal{I}$ we have
    $(\beta\tens\gamma)_- = 0$ hence $\beta_1 = \beta\tens\gamma$. More generally, $\beta
    \tens \alpha_k$ is irreducible and equivalent to $\beta_k$ for every $k$, so that the
    restriction of $\Te\tens\id$ to $p_\beta H \tens P_\gamma K$ is an isometry onto
    $P_{\beta,\gamma} K$: this is a trivial case of Proposition~\ref{prp:but_inj} and
    Remark~\ref{rque:ext_target}. For the same reason one has $(\beta_k \tens
    \bar\gamma_k)_+ \simeq \beta\tens (\alpha_k \tens \bar\gamma_k)_+$, which implies that
    $p_{\ss\jok+} (\Te \tens\id) = (\Te \tens\id) (\id\tens p_{\ss\jok+})$, and the
    similar relations for $p_{\ss+\jok}$. Moreover we also have $\Theta(\Te\tens\id) =
    (\Te\tens\id) (\id\tens\Theta)$ because $S\tens 1$ commutes to $\Theta$.
  \end{enumerate}
\end{rque}

\begin{thm} \label{thm:proj2}
  Let $S$ be a Woronowicz \Cst algebra and $p_1$ a central projection of $\hat S$ such
  that $Up_1U = p_1$ and $p_0 p_1 = 0$. Assume that the classical Cayley graph $\GG$ of
  $(S, p_1)$ is a directional tree.
  \begin{enumerate}
  \item \label{item:lim_inj} The maps $R_k$ are injective. As a result, the space
    $K_\infty$ is infinite-dimensional whenever $S$ is not co-commutative.
  \item \label{item:proj_closed} If we have $M_\gamma \neq 2$ for all $\gamma\in\Dir$,
    then $p_{\ss++}K_g = \Ker Rs$. Otherwise $p_{\ss++}K_g$ is a strict, dense subspace of
    $\Ker Rs$.
  \end{enumerate}
\end{thm}

\begin{dem} 
  \ref{item:lim_inj}. Thanks to the preceding
  Remark~\ref{rque:subtrees}.\ref{item:gen_subtree}, it is enough to study the restrictions of
  the considered objects to the subspaces $P_\gamma K$ with $\gamma\in \Dir$.  Let $l\in\NN$,
  we can suppose that $l\in\inter{1}{k}$, and in particular that $\dim\gamma > 1$: otherwise
  $p_{\ss+-} (p_k\tens p_1) q_l P_\gamma = 0$ hence $q_l P_\gamma K$ doesn't meet the
  definition set of $R_k$.  Because the subspaces $(p_k\tens\id) p_{\ss+-} q_l P_\gamma$ are
  irreducible, and by definition of the norm of $K_{\infty}$, we have
  \begin{displaymath}
    ||R_k q_l P_\gamma|| = \lim || r^{k+i} (p_k\tens\id) q_l P_\gamma || = 
    \prod_{i=0}^\infty || p_{\ss+-}\Theta p_{\ss+-} (p_{k+i}\tens\id) q_l P_\gamma || \text{.}
  \end{displaymath}
  
  To prove that this infinite product is non-zero we use the quantitative result of
  Lemma~\ref{lem:lecture_theta}. Recall from Lemma~\ref{lem:conjug_morph} that the
  sequence $(m_k)$ satisfies the induction equation $m_{i-1} - m_1m_i + m_{i+1} = 0$, so
  that we can write $m_i = (a^{i+1} - a^{-i-1}) / (a-a^{-1})$ for some $a>1$ when $m_1>2$,
  and $m_i = i+1$ when $m_1 = 2$. It is now very easy to check that the following infinite
  sum is finite:
  \begin{displaymath}
    \log \prod_{i=k}^\infty || p_{\ss+-}\Theta p_{\ss+-} (p_i\tens\id) q_l P_\gamma ||
    = \frac12 \sum_{i=k}^\infty \log \left(1 - \frac {m_lm_{l-1}} {m_{i+1}m_i}\right) 
    > -\infty \text{.}
  \end{displaymath}
  Note that we have $||R_k q_l P_\gamma|| \to 1$ when $k\to\infty$, and in particular the
  norm $||q_l P_\gamma R_k^{-1}||$ is bounded with respect to $k$. We will need to know
  for the second point that it is even bounded with respect to $k$ and $l$, when $m_1 >
  2$. To see this, check that $m_l / m_i \leq a^{-(i-l)}$ when $l\leq i$ and conclude that
  \begin{displaymath}
    \forall k\geq 1,~ l\in\inter{1}{k}~~ \log ||R_k q_l P_\gamma|| \geq 
    \ts\frac 12 \sum_{i=1}^{\infty}\ds \log \left(1-a^{-2i}\right) \text{.}
  \end{displaymath}
  
  Now if there indeed exists a direction $\gamma\in\Dir$ with $\dim \gamma > 1$, the
  injectivity of $R_k$ implies that $\dim K_{\infty} > \dim (p_k\tens p_1) p_{\ss+-}
  P_\gamma K = \dim \alpha_{k-1} \dim \alpha_{k+1}$, which tends to infinity with $k$
  according to Lemma~\ref{lem:strict}.

  \bigskip
  
  \ref{item:proj_closed}. We will use the decomposition given by the $q_l P_\gamma$ to
  study the expression of $p_{\ss++} K_g$ obtained in Proposition~\ref{prp:RS}, and in
  particular the operator $Rs p_{\geq k}$ restricted to $\Ker Rs$. If $l=0$ we have
  $p_{\ss++} P_\gamma K_g = P_\gamma (\Ker Rs) = P_\gamma K_{\ss++}$ since we are
  considering a classical graph. Now we assume that $l \geq 1$. In particular the map $s :
  (p_k\tens\id) q_l P_\gamma K_{\ss++} \to (p_{k+1}\tens\id) q_l P_\gamma K_{\ss+-}$ is
  bijective for any $k\geq l$ according to Lemma~\ref{lem:lecture_theta}, and hence
  $R_{k+1}s : (p_k\tens\id) q_l P_\gamma K_{\ss++} \to K_\infty$ is injective. Therefore
  it is possible to unitarily identify all the subspaces $(p_k\tens\id) q_l P_\gamma
  K_{\ss++}$ to their common image $G_l \subset K_\infty$ in order to have $R_{k+1}s q_l
  P_\gamma = \lambda_{k,l} \id_{G_l}$ with
  \begin{displaymath}
    \lambda_{k,l} = ||R_{k+1} s q_l P_\gamma|| = ||R_{k+1} q_l P_\gamma|| 
    \sqrt{ \frac {m_lm_{l-1}} {m_{k+1}m_k}} \text{.}
  \end{displaymath}
  In particular the operator $((p_k\tens\id)\zeta)_k \mapsto (Rsp_{\geq k}\zeta)_k$ from
  $q_l P_\gamma K_{\ss++}$ to $G_l^\NN$ identifies then with the augmentation by $G_l$ of
  the matrix $\Lambda_l = (\lambda_{j,l}\delta_{j\geq i\geq l})_{i,j}$.

  \bigskip

  \ref{item:proj_closed}a. We start with the case $m_1 = M_\gamma > 2$, which is particularly
  simple. As a matter of fact, $\Lambda_l$ is then bounded, even as an operator from
  $\ell^\infty(\NN)$ to $\ell^2(\NN)$, and uniformly with respect to $l$: we have
  \begin{eqnarray*}
    \sum_{i\geq l} \left( \ts\sum_j\ds |\lambda_{j,l}\delta_{j\geq i}| \right)^2 &=& 
    \sum_{i\geq l} \left( \ts\sum_{j\geq i}\ds  ||R_{j+1} q_l P_\gamma|| 
      \sqrt{ \ts\frac {m_lm_{l-1}} {m_{j+1}m_j} \ds} \right)^2 \\
    &\leq& \sum_{i\geq l} \left( \ts\sum_{j\geq i}\ds  a^{-(j+1-l)} \right)^2
      = \ts\frac {a^2} {(a^2 - 1) (a-1)^2}\ds \text{,}
  \end{eqnarray*}
  using the same estimate for $m_l/m_j$ as in the first point. As a result, for any vector
  $\zeta \in P_\gamma K_{\ss++}$ the sequence $(Rsp_{\geq k} \zeta)_k$ is square-summable. The
  operators $R_k^{-1}$ being uniformly bounded in our case, the sequence $(R_k^{-1} Rsp_{\geq
    k} \zeta)_k$ is also square-summable. Therefore the condition of Proposition~\ref{prp:RS}
  is satisfied by any vector in $P_\gamma(\Ker Rs)$.

  \bigskip
  
  \ref{item:proj_closed}b. Now we address the case $m_1 = 2$. Let $\varepsilon > 0$, there
  exists $I\geq l$ such that $||R_i q_l P_\gamma|| \geq 1-\varepsilon$ for every $i \geq
  I$. We have then the following inequalities:
  \begin{displaymath}
    \begin{array}{rrrl}
      & \delta_{j\geq i\geq I} (1-\varepsilon) \sqrt{\ts \frac {(l+1)l} {(j+2)(j+1)}} \leq
      &\lambda_{j,l} \delta_{j\geq i\geq l}& \leq 
      \delta_{j\geq i} \sqrt{\ts \frac {(l+1)l} {(j+2)(j+1)}} \\[1.5ex]
      \Longrightarrow & \delta_{j\geq i\geq I} (1-\varepsilon) {\ts\frac{l}{i+j+2}} \leq
      &\lambda_{j,l} \delta_{j\geq i\geq l}& \leq {\ts\frac{2l+2}{i+j+2}} \\[1.5ex]
      \Longrightarrow & l(1-\varepsilon)~ [ \delta_{i,j\geq I}\, \mu_{i+1,j+1}] \leq 
      &\Lambda_l + \Lambda_l^*& \\[1ex]
      & \text{and} & \Lambda_l &\leq (2l+2)~ [ \mu_{i,j} ] \text{,}
    \end{array}
  \end{displaymath}
  where we put $\mu_{i,j} = (i+j+1)^{-1}$. The last two inequalities are understood in the
  coefficientwise meaning, but it is well known that this implies norm inequalities,
  because all the coefficients are non-negative. Hence we have
  \begin{displaymath}
    {\ts\frac{l(1-\varepsilon)}2}~ || [\delta_{i,j\geq I}\, \mu_{i+1,j+1}] || \leq {\ts\frac12} 
    ||\Lambda_l + \Lambda_l^*|| \leq ||\Lambda_l|| \leq (2l+2)~ || [\mu_{i,j}] || \text{.}
  \end{displaymath}
  Now we have in the left-hand (resp. right-hand) side a compact perturbation of (resp.
  exactly) of the Hilbert matrix $M = [\mu_{i,j}]$, which is known from the theory of
  Hankel operators to have a norm and an essential norm both equal to $\pi/2$
  (cf \cite{nikolskij:toeplitz}, th.~5.3.1). Hence we obtain, letting furthermore
  $\varepsilon$ go to zero, the estimate $l\pi/4 \leq ||\Lambda_l|| \leq (l+1)\pi$.
  
  From this we conclude that every vector of $q_l P_\gamma K_{\ss++}$ satisfies the condition
  of Proposition~\ref{prp:RS} --- recall that the operators $q_l P_\gamma R_k^{-1}$ are
  uniformly bounded with respect to $k$. As a result, $p_{\ss++}K_g$ is dense in $\dirsum q_l
  P_\gamma (\Ker Rs) = \Ker Rs$.  However, $p_{\ss++}K_g$ is not equal to $\Ker Rs$. As a
  matter of fact, the lower estimate we have obtained proves that there exist vectors $\zeta_l
  \in q_l P_\gamma K_{\ss++}$ such that $||\zeta_l|| = 1/l$ and $||(Rsp_{\geq k} \zeta_l)_k||
  \geq \pi/4$. Moreover one can assume that $Rs(\zeta_l) = 0$: this only corresponds to
  composing $\Lambda_l$ on the right by a co-rank $1$ projection, which is a compact
  perturbation. One has then $\zeta = \sum_l \zeta_l \in P_\gamma(\Ker Rs)$, but $(Rsp_{\geq k}
  \zeta)_k$ is not square-summable.
\end{dem}

\section{Edges at infinity: the action}
\label{section:infinity_action}

In the previous section, the interest of the Hilbert space $K_{\infty}$ mainly lay in its
relation with the Hilbert space of geometric edges, via the projection $p_{\ss++}$. The
aim of this section is to endow $K_\infty$ with a representation of $S_\red$, which will
turn it into an interesting geometric object on its own. On the way, we will be led to
study certain aspects of the regular representation $S_\red \subset L(H)$ which can be of
independent use: see Lemma~\ref{lem:reg_asc_commut} and the remarks after it.

A first step however will be to notice that $K_\infty$ can easily be equipped with a
representation of $\hat S$, namely the inductive limit $\hat\pi_\infty$ of the representation
$\hat\pi_4 \hat\delta^3$. As a matter of fact the image of $\hat\pi_4 \hat\delta^3$ commutes to
$p_{\ss+-}$, $p_k\tens\id$ and $r$. Recall from the preceding section that the decomposition of
$(p_k\tens\id)K_{\ss+-}$ into irreducible subspaces with respect to $\hat\pi_4 \hat\delta^3$
are given by the projections $P_{\beta,\gamma} q_l$.  The subspace $P_{\beta,\gamma} q_l
(p_k\tens\id) K_{\ss+-}$ is non-zero \iff $\dim\gamma > 1$ and $|\beta|+1 \leq l\leq k$ and is
then equivalent to $\beta \tens \gamma \bar\gamma \cdots \gamma_{2l-2|\beta|} \tens \bar\beta$.
As a result $\hat\pi_\infty (p_\alpha) K_\infty$ is irreducible if $\alpha\neq 1_\Cat$ and
$\alpha\subset \delta\tens \bar\delta$ for some $\delta\in \Irr\Cat$, and vanishes else.

\begin{lem} \label{lem:reg_asc_commut}
  Let $p_{\ss\jok+}$ be the left ascending projection of a quantum Cayley tree. Choose
  $\gamma\in\Dir$ with $\dim\gamma > 1$ let $m_k = M_{\alpha_k}$ be the corresponding
  sequence of quantum dimensions. Let $a \in S_\red$ be a coefficient of $\gamma$. Then the
  commutator $[a\tens 1, p_{\ss\jok+}]$ vanishes on $(1-P_{\bar\gamma})K$, and there
  exists a real number $C_a>0$ such that
  \begin{displaymath}
    \forall~ k\in\NN~~ ||[a\tens 1, p_{\ss\jok+}] (p_k\tens\id)|| \leq {C_a} {m_k}^{-1} \text{.}
  \end{displaymath} 
\end{lem}

\begin{dem}
  For this proof we can of course assume that $k$ is greater than $2$. It is enough to
  study $p_{\ss\jok+} [a\tens 1, p_{\ss\jok+}] (p_k\tens\id)$ because $[a\tens 1,
  p_{\ss\jok+}] = [a\tens 1, p_{\ss\jok+}] p_{\ss\jok+}$ $- p_{\ss\jok+} [a\tens 1,
  p_{\ss\jok+}]$. We will use the ``extended target operator'' $\Te : p_\gamma H\tens H
  \to H$ given by the product of $S$. Denoting by $\tilde a$ the map $(\CC\to H, 1\mapsto
  \Lambda_h(a))$ we have
  \begin{eqnarray*}
    p_{\ss\jok+} [a\tens 1, p_{\ss\jok+}] &=& p_{\ss\jok+} (\Te\tens\id) 
    (\tilde a\tens\id_K)  p_{\ss\jok+} - p_{\ss\jok+} (\Te\tens\id)(\tilde a\tens\id_K)  \\ 
    &=& (\Te\tens\id)~ (\hat\delta\tens\id)(p_{\ss\jok+})~
    (\id\tens p_{\ss\jok+} - 1)~ (\tilde a\tens\id_K) \text{.}
  \end{eqnarray*}
  Hence it is enough to show that $||P_1(1-P_2)|| \leq m_k^{-1}$, where $P_1$ and
  $P_2$ are the respective restrictions to $p_\gamma H\tens p_k H\tens p_1 H$ of
  $(\hat\delta\tens\id) (p_{\ss\jok+})$ and $(\id\tens p_{\ss\jok+})$. These projections
  act through the left representation of $\hat S^{\tensexp 3}$ on $H^{\tensexp 3}$, so
  that it suffices to look at their action on $L = H_\gamma\tens H_\alpha\tens
  H_{\gamma'}$, with $\gamma' \in\Dir$ and $|\alpha| = k$.
  
  Let $H_{(\gamma\tensexp\alpha)_+}$ and $H_{(\gamma\tensexp\alpha)_-}$ be the irreducible
  subspaces of $H_\gamma \tens H_\alpha$, the latter being possibly vanishing. We let
  $p_{\ss\jok+} \in M(\hat S\tens\hat S)$ act on any representation space of $\hat
  S\tens\hat S$. The image of $P_2$ is then $H_\gamma\tens p_{\ss\jok+} (H_\alpha \tens
  H_{\gamma'})$, whereas the image of $P_1$ is the sum of $L_+ = p_{\ss\jok+}
  (H_{(\gamma\tensexp\alpha)_+} \tens H_{\gamma'})$ and $L_- = p_{\ss\jok+}
  (H_{(\gamma\tensexp\alpha)_-} \tens H_{\gamma'})$. Let us first consider the case when
  $\alpha$ is not one of the form $\delta\bar\delta \cdots \delta_k$ for any
  $\delta\in\Dir$. Notice that we are automatically in this case when $\dim\gamma = 1$,
  because we restricted ourselves to the values $k\geq 2$. The corepresentation $\alpha$
  can then be written as an irreducible tensor product $\alpha_1\tens\alpha_2$, so that
  one has
  \begin{eqnarray*}
    L_+ &=& p_{\ss\jok+} ((H_{(\gamma\tensexp\alpha_1)_+} \tens H_{\alpha_2}) 
    \tens H_{\gamma'}) \\ &=& (\id\tens p_{\ss\jok+}) (H_{(\gamma\tensexp\alpha_1)_+} 
    \tens (H_{\alpha_2} \tens H_{\gamma'})) \\ &\subset& 
    H_\gamma\tens H_{\alpha_1} \tens~ p_{\ss\jok+} (H_{\alpha_2} \tens H_{\gamma'}) =
    H_\gamma\tens p_{\ss\jok+} (H_\alpha \tens H_{\gamma'}) \text{,}
  \end{eqnarray*}
  and similarly $L_- \subset \Im P_2$. In this case we therefore have $(1-P_2) P_1 = 0$.
  One can check in the same way that it is also the case when the geodesic from $1_\Cat$
  to $\alpha$ does not start in the direction $\bar\gamma$ or does not end with the
  direction $\bar\gamma'$.
  
  Therefore it remains to consider the situation when $\gamma$ is the generator of some
  copy of $A_o(Q)$ or $A_u(Q)$ in $S$, and $\alpha = \bar\gamma \gamma \cdots
  \bar\gamma_k$, $\gamma' = \gamma_k$. In other words we have $H_\gamma \tens
  H_\alpha\tens H_{\gamma'} = H_{1,k,1}$ with the notation of Lemma~\ref{lem:angles_2}.
  Let us notice first that $L_+$ is the unique irreducible subspace of $H_{1,k,1}$ which
  is at distance $k+2$ from the origin $1_\Cat$, and is therefore included in $(\id\tens
  p_{\ss\jok+})(H_{1,k,1})$. Hence it suffices to consider the restriction of $P_1$ and
  $P_2$ to the copies of $H_k$ in $H_{1,k,1}$. We are then exactly in the situation of
  Lemma~\ref{lem:angles_2}, with $k'=1$, $G_1 = \Im P_1$ and $G_2 = \Im P_2$. Because we
  are looking now at morphisms between irreducible subspaces, we can finally use the lemma
  to write
  \begin{displaymath}
    ||(1 - P_2) P_1||^2 = 1 - ||P_2 P_1||^2 = {m_k^{-2}} \text{.}
  \end{displaymath}
\end{dem}

\begin{rque}{Remarks} \label{rque:reg_commut} 
  \begin{enumerate}
  \item \label{item:min_quant_dim} Let $(m_k)_k$, $(m'_k)_k$ be two sequences of quantum
    dimensions associated to two directions $\gamma$, $\gamma' \in \Dir$. If $m'_1 \geq
    m_1$, it is easy to check by induction, using point~\ref{item:conjug_induc_dim} of
    Lemma~\ref{lem:conjug_morph}, that $m'_{k+1}/m'_k \geq m_{k+1}/m_k$:
    \begin{eqnarray*}
      && m'_{k+1}m_k - m_{k+1}m'_k = \\ && \makebox[1.5cm]{} 
      = (m'_1-m_1)m'_km_k + (m'_km_{k-1} - m_km'_{k-1}) \geq 0 \text{.} 
    \end{eqnarray*}
    In particular we have $m'_k \geq m_k$ for all $k$. If $m_1$ is minimal (resp. maximal)
    amongst the $M_\gamma$ with $\gamma\in\Dir$ and $\dim\gamma > 1$, we will call
    $(m_k)_k$ the minimal (resp. maximal) sequence of quantum dimensions for $(S,p_1)$.
  \item \label{item:reg_asc_commut_gen} It is clear from the proof of the lemma that
    $[p_{\ss\jok+}, a\tens 1] (p_k\tens\id)$ vanishes as soon as $k\geq 2$ if $a$ is a
    coefficient of some $\gamma\in\Dir$ with $\dim\gamma = 1$. Let us prove now that the
    result of the lemma holds in fact for any $a \in \mathcal{S} \subset S_\red$ if one
    uses the minimal sequence of quantum dimensions to state it. To see this, assume that
    $a$ satisfies the inequalities of the lemma and let $u$ be a coefficient of a
    corepresentation $\gamma\in\Dir$. Because the algebra $\Ss$ is spanned by such
    coefficients, it is enough to prove that $au$ also satisfies the same inequalities for
    some other constant $C_{au}$. We remark that $(u\tens 1) (p_k\tens\id) K$ is included
    in $(p_{k-1}\tens\id)K + (p_{k+1}\tens\id)K$, so that one can write, using the
    inequalities $m_k \leq m_{k+1} \leq m_1m_k$:
    \begin{eqnarray*} 
      && \makebox[-1cm]{} || [au\tens 1, p_{\ss\jok+}] (p_k\tens\id) || \leq \\
      && \leq ||(a\tens 1) [u\tens 1, p_{\ss\jok+}] (p_k\tens\id)|| + 
      || [a\tens 1, p_{\ss\jok+}] (u\tens 1) (p_k\tens\id)|| \\
      && \leq ||a|| C_u m_k^{-1} + ||u|| C_a (m_{k-1}^{-1} + m_{k+1}^{-1}) \\
      && \leq (||a|| C_u + (m_1 + 1) ||u|| C_a)~ m_k^{-1} \text{.} 
    \end{eqnarray*}
  \item \label{item:reg_asc_commut_ext} The lemma also admits the following
    generalization.  If we put $\Pe_{\ss\jok+} = \sum (p_k\tens p_{k'})
    \hat\delta(p_{k+k'})$ as in Remark~\ref{rque:ext_target}, we have for any coefficient
    $a$ of any $\gamma\in\Dir$ and for any $k$, $k'\in\NN^*$:
    \begin{equation} \label{eq:reg_asc_commut_ext}
      || [a\tens 1, \Pe_{\ss\jok+}] (p_k\tens p_{k'})|| \leq 
      C_a \sqrt{\ts\frac {m_{k'-1}}{m_{k+k'-1} m_k}} \text{,}
    \end{equation}
    where $(m_k)$ is the sequence of quantum dimensions associated with $\gamma$.
    Moreover $[a\tens 1, \Pe_{\ss\jok+}] (p_\alpha \tens p_\beta)$ can only be non-zero
    when $\alpha = \bar\gamma\gamma \cdots \bar\gamma_k$ and $\beta = \gamma_k\cdots
    \gamma_{k+k''-1}\tens \beta'$ with $k''\geq 1$ --- we have then $k' = k'' + |\beta'|$.
    Notice that in this case $M_\alpha = m_{|\alpha|}$ and the subobject $\delta \subset
    \alpha\tens \beta$ with maximal length is $\bar\gamma \cdots \bar\gamma_{k+k''}\tens
    \beta'$, so that $M_\beta / M_\delta$ equals $m_{k''} / m_{k+k''}$, which is less than
    $m_{|\beta|} / m_{|\delta|}$. We will use these facts in the proof of
    Theorem~\ref{thm:AO}.
    
    To prove~(\ref{eq:reg_asc_commut_ext}), one considers like in the proof of the lemma
    intertwining projections in $H_\gamma\tens H_\alpha\tens H_\beta$: the complete
    statement of Lemma~\ref{lem:angles_2} gives then the result with $m_{k''-1}/
    m_{k+k''-1}$. But this quotient is less than $m_{k'-1}/ m_{k+k'-1}$, because the
    sequence $(m_{k'-1}/ m_{k+k'-1})_{k'}$ is non-decreasing for every $k$: compare $m_{k'-1} m_{k+k'}$ and $m_{k+k'-1} m_{k'}$ by considering the irreducible
    decompositions of $H_{k'-1,k+k'}$ and $H_{k+k'-1,k'}$ relative to the appropriate
    $A_o(Q)$ or $A_u(Q)$.
  \item \label{item:reg_finite_propag} Using the same starting point as in the proof of
    Lemma~\ref{lem:reg_asc_commut}, we can prove the following result: if $a\in\Ss\subset
    S_\red$ is a coefficient of the corepresentation $\alpha \in \Irr\Cat$, and for any
    $\beta\in\Irr\Cat$, we have $a p_\beta H \subset \sum \{ p_\delta H ~|~ \delta \subset
    \alpha\tens\beta \}$. As a matter of fact one can write, using the notation of the
    proof, $p_\delta a p_\beta = p_\delta \Te (\tilde a \tens p_\beta) = \Te
    \hat\delta(p_\delta) (\tilde a \tens p_\beta)$. But $\tilde a$ lies in $p_\alpha H$ by
    assumption, hence the considered product vanishes if $\delta \not\subset
    \alpha\tens\beta$.  Similarly, one can check that $(a\tens 1)\hat\pi_4
    \hat\delta^3(p_\beta) K$ is included in the sum of the $\hat\pi_4
    \hat\delta^3(p_\delta)K$ with $\delta \subset \alpha\tens\beta\tens\alpha$.  These
    ``propagation properties'' will in particular be used in relation with the following
    elementary fact: if $H = \dirsum pH = \dirsum qH$ are orthogonal decompositions of
    $H$, and if $f\in L(H)$ is an operator such that $\mathrm{Card}~ \{q ~|~ pfq \neq 0\}
    \leq N$ for all $p$, one has $||f|| \leq \sqrt N~ \sup ||fq||$.
  \end{enumerate} 
\end{rque}

\begin{thm} \label{thm:infinite_reg}
  Let $S$ be a Woronowicz \Cst algebra and $p_1$ a central projection of $\hat S$ such
  that $Up_1U = p_1$ and $p_0 p_1 = 0$. Assume that the classical Cayley graph $\GG$ of
  $(S, p_1)$ is a directional tree. Let us denote by $\varphi_{\ss+-}(a)$ the operator $p_{\ss+-}
  (a\tens 1) p_{\ss+-}$, for any $a\in S_\red$.
  \begin{enumerate}
  \item Let $\zeta\in (p_k\tens p_1)K_{\ss+-}$ and $a\in\Ss \subset S_\red$. The sequence
    $(R \varphi_{\ss+-}(a) r^n \zeta)_n$ converges in $K_\infty$ to a vector which only depends 
    on $R\zeta$ and which we denote by $\pi_\infty(a) (R\zeta)$.
  \item This defines a $*$-algebra morphism $\pi_\infty : \Ss \to L(K_\infty)$ which extends by
    continuity to $S_\red$.
  \end{enumerate}
\end{thm}

\begin{dem}
  Let $a$ be an element of $\Ss\subset S_\red$. There exists an integer $p$ such that $a$
  can be expressed as a sum of coefficients of corepresentations $\beta\in\Irr\Cat$ with
  $|\beta| \leq p$. We will use in this proof the finite propagation properties of $a$,
  see Remark~\ref{rque:reg_commut}.\ref{item:reg_finite_propag}, with respect to two
  decompositions of $K$. The first one is simply given by the projections $(p_k\tens\id)$,
  but the second one is a little bit more subtle. Using the notation of
  Remark~\ref{rque:subtrees}.\ref{item:gen_subtree}, for $k_0\in\NN$ and $l\in\NN^*$ we
  denote by $Q_{k_0, l}$ the sum of the projections $P_{\beta,\gamma} q_{k_0 + l}
  p_{\ss+-}$ with $(\beta,\gamma) \in \mathcal{I}$ and $|\beta| = k_0$. In other words,
  the $\hat S$-subspace $(p_{k+k_0}\tens\id) Q_{k_0,l} K$ is the sum over $(\beta,\gamma)
  \subset \mathcal{I}$, $|\beta| = k_0$, of the irreducible subspaces $\beta\tens
  \alpha_{2l} \tens\bar\beta \subset (p_{\beta\tensexp\alpha_k} \tens p_{\bar\gamma_k})
  K_{\ss+-}$, where $\alpha_k = \gamma \cdots \gamma_k$ as usual. In particular we have $R 
  Q_{k_0,l} = Q^\infty_{k_0,l} R$, if $Q^\infty_{k_0,l}$ is the sum of the projections
  $\hat\pi_\infty (p_{\beta \tensexp\alpha_{2l} \tensexp\bar\beta})$.

  \bigskip
 
  We first want to bound from above the norm of the commutator $\Cc_a = [\varphi_{\ss+-}(a),
  r]$ on each subspace $Q_{k_0,l} K$. Because $S\tens 1$ commutes to $\Theta$, and using
  Proposition~\ref{prp:orient}, we see that the operator $p_{\ss+-}(a\tens 1)p_{\ss\jok-}
  \Theta p_{\ss+-}$ equals $p_{\ss+-}\Theta p_{\ss+\jok} (a\tens 1) p_{\ss+-}$. We subtract and
  add this quantity from $\Cc_a$ and force the apparition of the commutators of
  Lemma~\ref{lem:reg_asc_commut}:
  \begin{eqnarray} \nonumber
     && \makebox[-2cm]{} \Cc_a = p_{\ss+-} (a\tens 1) p_{\ss+-} \Theta p_{\ss+-} - 
     p_{\ss+-} \Theta p_{\ss+-} (a\tens 1) p_{\ss+-} = \\ \nonumber 
     && = -p_{\ss+-} (a\tens 1) p_{\ss--} \Theta p_{\ss+-}
     + p_{\ss+-}\Theta p_{\ss++} (a\tens 1) p_{\ss+-} \\ \label{eq:reg_shift_commut}
     && = p_{\ss+-}[a\tens 1, p_{\ss+\jok}] p_{\ss--} \Theta p_{\ss+-}
     - p_{\ss+-}\Theta p_{\ss++} [a\tens 1, p_{\ss\jok +}] p_{\ss+-} \text{.}
  \end{eqnarray}
  
  Thanks to Remark~\ref{rque:subtrees}.\ref{item:gen_subtree}, the norm of $p_{\ss--}
  \Theta p_{\ss+-}$ on $(p_k\tens\id) Q_{k_0, l} K$ is the same as the one on
  $(p_{k-k_0}\tens\id) Q_{0,l} K$, which is given by Lemma~\ref{lem:lecture_theta}:
  \begin{equation} \label{eq:reg_shift_theta}
    || p_{\ss--} \Theta p_{\ss+-} (p_k\tens\id) Q_{k_0, l} ||^2 \leq
    \frac {M_l M_{l-1}} {m_{k-k_0+1} m_{k-k_0}} \text{,}
  \end{equation}
  where $(m_k)_k$ and $(M_k)_k$ are the minimal and maximal sequences of quantum
  dimensions from Remark~\ref{rque:reg_commut}.\ref{item:min_quant_dim}.  We proceed in
  the same way for the second term of~(\ref{eq:reg_shift_commut}), but this time we have
  to consider the restriction of $p_{\ss+-}\Theta p_{\ss++}$ to the subspaces
  $(p_{k'}\tens\id) Q_{k'_0,l'} K$ that meet the image of $(a\tens 1) (p_k\tens\id)
  Q_{k_0,l}$.  Remark~\ref{rque:reg_commut}.\ref{item:reg_finite_propag} provides control
  over the set of indices $(k', k'_0, l')$ to be considered, and the fact that quantum
  dimensions are increasing with the distance to the origin shows that the greatest value
  of the quantity~(\ref{eq:reg_shift_theta}) is obtained when $(k', k'_0, l') = (k-p,
  k'_0+p, l+p)$.  Putting this together with the estimate of
  Lemma~\ref{lem:reg_asc_commut} we get
  \begin{displaymath}
    ||\Cc_a (p_k\tens\id) Q_{k_0,l}|| \leq \ts 
    2~ \frac {C_a} {m_k}~ \sqrt{\frac{M_{l+p} M_{l+p-1}} {m_{k-k_0-2p+1}m_{k-k_0-2p}}}
    \ds \leq \frac {C_{a,k_0,l}} {m_k^2} \text{.}
  \end{displaymath}
  Notice that we have used the inequality $m_{k-i} \geq m_k m_1^{-i}$ and introduced a new
  constant $C_{a,k_0,l}$ to obtain the estimate order $m_k^{-2}$.

  \bigskip
  
  Now we consider a vector $\zeta \in (p_k\tens\id) Q_{k_0,l} K$, for fixed integers $k_0$
  and $l$. To prove that the sequence $(R \varphi_{\ss+-}(a) r^n \zeta)$ converges, it is
  enough to study the series $(\sum R \varphi_{\ss+-}(a) r^{n+1} \zeta - R
  \varphi_{\ss+-}(a) r^n \zeta)$, which can be written as $(\sum R \Cc_a r^n \zeta)$.
  Because the vector $\Cc_a r^n \zeta = [\varphi_{\ss+-}(a), r] r^n \zeta$ belongs to the
  direct sum of the subspaces $(p_{k+n+i+1}\tens\id) K_{\ss+-}$ with $i\in \inter{-p}{p}$,
  we have
  \begin{displaymath}
    || R \Cc_a r^n \zeta || \leq (2p+1) || \Cc_a r^n \zeta || \leq
    \frac {2p+1} {m_{k+n}^2} ~ C_{a,k_0,l} ||\zeta|| \text{.}
  \end{displaymath}
  Now we have $m_{k+n} \geq k+n+1$, hence the series $(\sum_n m_{k+n}^{-2})$ is convergent
  and the sequence $(R\varphi_{\ss+-}(a) r^n \zeta)_n$ indeed converges in $K_\infty$. If
  $R\zeta = R\zeta'$ with $\zeta' \in (p_{k'}\tens\id) K_{\ss+-}$ and $k'\geq k$, we have
  $\zeta' = r^{k'-k} \zeta$ by injectivity of $R_{k'}$, hence the associated sequences are
  equal up to an index shift.
  
  We moreover get an estimate on the norm of $||\pi_\infty(a) RQ_{k_0,l}||$: denoting by
  $(\rho_i)_i$ the sequence of remainders of the series $(\sum m_i^{-2})$, we have
  \begin{eqnarray} \label{eq:infinite_reg_conv}
    & || \pi_\infty(a)(R\zeta) - R \varphi_{\ss+-}(a) \zeta || \leq
    C_{a,k_0,l}~ (2p+1) \rho_{k} ||\zeta|| \text{,} \\ \nonumber \text{hence}
    & || \pi_\infty(a)(R\zeta)|| \leq (2p+1) (||a|| + C_{a,k_0,l}\rho_{k}) ||\zeta|| \text{.}
  \end{eqnarray}
  If we let $k$ go to infinity without changing $R\zeta$, the norm of $\zeta$ converges to
  $||R\zeta||$ and we get the upper bound $||\pi_\infty(a) Q^\infty_{k_0,l}|| \leq (2p+1)
  ||a||$.  We finally use Remark~\ref{rque:reg_commut}.\ref{item:reg_finite_propag} to
  notice that $\varphi_{\ss+-} (a) Q_{k_0,l} K$ is included in the sum of the $(2p+1)^2$
  subspaces $\{Q_{k_0+i_0,l+j}K ~|~ i_0 \text{~and~} i_0+j \in \inter{-p}{p}\}$. As a
  result the same property of ``finite propagation'' is true for $\pi_\infty(a)$ in the
  decomposition $K_\infty = \dirsum Q^\infty_{k_0,l} K_\infty$ and we obtain
  the inequality $||\pi_\infty(a)|| \leq (2p+1)^2 ||a||$.

  \bigskip
  
  Let $a$, $a' \in \Ss \subset S_\red$ and $R_k \zeta \in K_\infty$. By
  Remark~\ref{rque:reg_commut}.\ref{item:reg_asc_commut_gen}, the norm $||(\varphi_{\ss+-}
  (a) \varphi_{\ss+-}(a') - \varphi_{\ss+-}(aa')) r^n\zeta||$ tends to zero as $n$ goes to
  infinity. By definition of $\pi_\infty$, the norm $||(\pi_\infty(b)R -
  R\varphi_{\ss+-}(b)) r^n \zeta||$, with $b = a$, $a'$ or $aa'$, also tends to zero. As a
  result, we see that $||(\pi_\infty(a)\pi_\infty(a') - \pi_\infty(aa')) R r^n \zeta||$
  converges to zero with respect to $n$. But this quantity does not depend on $n$, and
  hence we have proved that $\pi_\infty$ is a morphism of algebras.
  
  In particular, it is enough to prove the identity $\pi_\infty(a^*) = \pi_\infty(a)^*$
  for the coefficients $a$ of any $\gamma\in\Dir$. We have then $\varphi_{\ss+-}(a)
  (p_{k+n}\tens\id) K \subset (p_{k+n+1}\tens\id)K + (p_{k+n-1}\tens\id)K$. Let
  $R\zeta$, $R\xi \in K_\infty$, we can assume that $\zeta$ and $\xi$ both lie in some
  $(p_k\tens\id) K_{\ss+-}$ and we write
  \begin{eqnarray*}
    (\pi_\infty(a^*)R\zeta | R\xi) &=& \lim~ (R\varphi_{\ss+-}(a)^* r^n\zeta | R\xi) \\
    &=& \lim~ (\varphi_{\ss+-}(a)^* r^n\zeta | r^{n+1}\xi) +
    \lim~ (\varphi_{\ss+-}(a)^* r^n\zeta | r^{n-1}\xi) \\
    &=& \lim~ (r^{n-1} \zeta | \varphi_{\ss+-}(a) r^n\xi)
    + \lim~ (r^{n+1} \zeta | \varphi_{\ss+-}(a) r^{n}\xi) \\
    &=& \lim~ (R\zeta | R \varphi_{\ss+-}(a) r^n\xi) = (R\zeta | \pi_{\infty}(a) R\xi) \text{.}
  \end{eqnarray*} 
  
  Finally, let us notice that if the ``propagation length'' of $a\in\Ss$ is $p$, the one
  of $a^n$ is at most $np$, so that $||\pi_\infty(a)^n|| \leq (2np+1)^2 ||a||^n$ for any
  $n$. In particular when $||a||<1$ in $S_\red$ this proves that $(\sum
  ||\pi_\infty(a)^n||)$ converges, so that the spectral radius of $\pi_\infty(a)$ is less
  than or equal to $1$. If $a$ is moreover hermitian, so is $\pi_\infty(a)$ and we get
  $||\pi_\infty(a)|| \leq 1$. Hence $\pi_\infty : \Ss \to L(K_\infty)$ is continuous when
  $\Ss$ is equipped with the norm of $S_\red$.
\end{dem}

\begin{rque}{Remark} \label{rque:simplif_infinite_reg}
  In the case when $M_\gamma \neq 2$ for all $\gamma\in\Dir$, the proof of the theorem can
  be simplified. More precisely, it is enough to use in~(\ref{eq:reg_shift_commut}) the
  evident upper bound $1$ for the norms of $p_{\ss--}\Theta p_{\ss+-}$ and
  $p_{\ss+-}\Theta p_{\ss++}$ --- and in particular there is no need to introduce the
  projections $Q_{k_0,l}$ anymore. As a matter of fact, the inequality
  $||\Cc_a(p_k\tens\id)|| \leq 2C_a m_k^{-1}$ is sufficient for the rest of the proof
  because the series $(\sum m_k^{-1})$ is geometrically convergent in this case.
\end{rque}

\section{Applications}
\label{section:appl}

\subsection{Property AO}
\label{section:AO}

In this section we will denote by $\lambda$ and $\rho : S \to L(H)$ the left and right
regular representations of a Woronowicz \Cst algebra $(S,\delta)$, ie $\rho(x) =
U\lambda(x)U$. They commute and therefore define a representation $(\lambda,\rho)$ of
$S_\red \tensmax S_\red$ on $H$.  Besides, we will call $\lambda\tens\rho$ the natural
representation of $S_\red \tensmax S_\red$ on $H\tens H$, so that
$(\lambda\tens\rho)(S_\red \tensmax S_\red) = S_\red \tens S_\red$. Let $\pi : L(H) \to
L(H)/K(H)$ be the quotient map. We say that $(S,\delta)$ has Property AO, after Akemann
and Ostrand, if $\pi \rond (\lambda, \rho)$ factorizes through $S_\red \tens S_\red$.

When the antipode of $(S,\delta)$ is involutive, it is easy to see that $(\lambda, \rho)
\rond \delta$ contains the trivial representation $\varepsilon$. Hence in this case
$(\lambda, \rho)$ factorizes through $S_\red \tens S_\red$ \iff $(S, \delta)$ is
amenable. Consequently, Property AO is only interesting for non-amenable Kac-\Cst algebras
and can be seen as a restriction on their non-amenability.

Property AO was first introduced in \cite{AkeOstrand:tensorfree} to study the non-nuclear
\Cst algebra $S = C^*(\FF_2)$: it was used in this case to show that $S_\red \tensmax
S_\red /$ $S_\red \tens S_\red \simeq K(H)$. This result was generalized to reduced
\Cst algebras of ICC discrete groups in \cite{Skand:knucl}, where Property AO was also
used in conjunction with Property T of Kazhdan to produce non-$K$-nuclear \Cst algebras.
More recently, Property AO was used in \cite{Ozawa:solid} in conjunction with local
reflexivity to produce solid factors.

\bigskip

The aim of this section is to prove Property AO for the free quantum groups $(S,\delta)$
studied in this article. We will use the original method of \cite{AkeOstrand:tensorfree}:
the factorization of $\pi\rond (\lambda,\rho)$ arises from an isometry $F : H \to H\tens
H$ such that $F^* (\lambda\tens\rho)(x) F \equiv (\lambda,\rho)(x)$ $\mathrm{mod}$ $K(H)$,
for any $x \in S_\red\tensmax S_\red$.  In the case of $\FF_2$, the isometry $F$ is the
polar part of the closable operator which maps each characteristic function
$\un_\alpha\in H$ to the sum of the $\un_{\beta_1} \tens \un_{\beta_2}$ with
$\beta_1\beta_2 = \alpha$ and $|\beta_1| + |\beta_2| = |\alpha|$. In particular, the
adjoint of this operator coincides on $\Hh\tens\Hh$ with the natural extension
$\Te\Pe_{\ss\jok+}$ of $\T p_{\ss\jok+}$. In the quantum case, we also define $F$ from
this extension.

\begin{df} \label{df:isom_AO}
  Let $S$ be a Woronowicz \Cst algebra and $p_1$ a central projection of $\hat S$ such
  that $Up_1U = p_1$ and $p_0 p_1 = 0$. Assume that the classical Cayley graph $\GG$ of
  $(S, p_1)$ is a directional tree. Let $\Te : \Hh\tens \Hh \to H$ be the operator induced in
  the GNS construction by the multiplication of $S$, and $\Pe_{\ss\jok+} = \sum_{k,k' \in \NN}
  \hat\delta(p_{k+k'}) (p_k\tens p_{k'})$. We define the closed operator $F_0$ by $F_0^* =
  \Te\Pe_{\ss\jok+}$ and we denote by $F$ its polar part.
\end{df}

\begin{lem} \label{lem:norm_AO}
  We use the hypotheses and notation of Definition~\ref{df:isom_AO}. We suppose that
  $M_\gamma \neq 2$ for all $\gamma\in\Dir$. For every $k\in\NN^*$ we have then $||F_0
  p_k||^2 \geq k+1$, and there exists a constant $C>0$ such that
  \begin{displaymath}
    \forall~ \alpha, \alpha'\in \Irr\Cat~~~ \alpha' \subset \Dir\tens\alpha ~~\Longrightarrow~~ 
    \left|\, ||F_0 p_{\alpha'}||^2 - ||F_0 p_{\alpha}||^2 \,\right| \leq C \text{.}
  \end{displaymath}
\end{lem}

\begin{dem}
  We know from Proposition~\ref{prp:but_inj} and Remark~\ref{rque:ext_target} that $F_0$
  is a multiple of an isometry on each subspace $p_\alpha H$, the corresponding norm being
  given by
  \begin{equation} \label{eq:norm_AO}
    || F_0 p_\alpha ||^2 = \sum \ts \Big\{ \frac {M_{\beta_1} M_{\beta_2}}{M_\alpha} ~\Big|~ 
    \alpha\subset \beta_1\tens\beta_2,~ |\beta_1| + |\beta_2| = |\alpha| \Big\} \text{.}
  \end{equation}
  Each term of the sum is clearly greater than or equal to $1$, and there are $|\alpha| +
  1$ terms in the sum: one obtains the admissible pairs $(\beta_1,\beta_2)$ by following
  the geodesic from $1_\Cat$ to $\alpha$ until an arbitrary point $\beta_1$, and then
  using the remaining sequence of directions to go up from $1_\Cat$ to $\beta_2$. Recall
  from Lemma~\ref{lem:strict} that the conditions for a sequence of directions to define
  an ascending path are only local.
  
  To get the second estimate, let us consider an inclusion $\alpha' \subset
  \gamma\tens\alpha$ with $\gamma \in \Dir$. By exchanging $\alpha$ and $\alpha'$ if
  necessary, one can assume that $|\alpha'| > |\alpha|$. As a first step, we will assume
  that $\alpha = \bar\gamma\gamma \cdots \bar\gamma_i$ and $\alpha' =
  \gamma\bar\gamma\gamma \cdots \gamma_i$. We can moreover suppose then that $i\geq 2$,
  and hence $\dim\gamma > 1$. Let $(m_k)$ be the sequence of quantum dimensions associated
  to $\gamma$, by hypothesis we have $m_k \sim \frac{a^{k+1}}{a-a^{-1}}$ for some $a>1$.
  We write then
  \begin{displaymath}
    f_i := || F_0 p_\alpha ||^2 = \sum_{k+k' = i} \frac {m_k m_{k'}} {m_i}
     = \frac {a^i}{m_i} \sum_{k+k' = i} \frac {m_k}{a^k}~ \frac {m_{k'}} {a^{k'}} \text{,}
  \end{displaymath}
  and similarly $|| F_0 p_{\alpha'} ||^2 = f_{i+1}$. But we have by a variant of Cesaro's
  Lemma $f_i \sim \frac {a(i+1)}{a-a^{-1}}$, and in particular $(f_{i+1} - f_i)_i$ is
  bounded. We take $C$ to be a common bound for these sequences when $\gamma$ varies in
  $\Dir$.
  
  We address now the general case and express $\alpha$ as a tensor product
  $\bar\gamma\gamma \cdots \bar\gamma_i \tens \tilde\alpha$, where $\tilde\alpha$ does not
  start with $\bar\gamma$, and possibly $i=0$ or $\tilde\alpha = 1_\Cat$. We have then
  $M_\alpha = m_i M_{\tilde\alpha}$ and $M_{\alpha'} =m_{i+1} M_{\tilde\alpha}$. Let us
  first consider the terms in~(\ref{eq:norm_AO}) where $|\beta_1| > i$. One has then
  $M_{\beta_1} = m_i M_{\tilde\beta_1}$ for some $\tilde\beta_1$, hence $M_{\beta_1}
  M_{\beta_2} / M_{\alpha} = M_{\tilde\beta_1} M_{\beta_2} / M_{\tilde\alpha}$. If we
  consider similarly in the expression~(\ref{eq:norm_AO}) for $||F_0 p_{\alpha'}||^2$ the
  terms where $|\beta_1| > i+1$, we see that $M_{\beta_1} = m_{i+1}M_{\tilde\beta_1}$ and
  $M_{\beta_1} M_{\beta_2} / M_{\alpha'} = M_{\tilde\beta_1} M_{\beta_2} /
  M_{\tilde\alpha}$. Hence all these terms can be simplified from the difference $||F_0
  p_{\alpha'}||^2 - ||F_0 p_{\alpha}||^2$.  We proceed symmetrically with the terms
  of~(\ref{eq:norm_AO}) where $|\beta_1| \leq i$ (resp.  $i+1$) : this time $\beta_2$ can
  be expressed as an irreducible tensor product $\tilde\beta_2 \tens \tilde\alpha$, the
  factors $M_{\tilde\alpha}$ disappear from the quotient $M_{\beta_1}M_{\beta_2}/M_\alpha$
  (resp. $M_{\alpha'}$) and one recognizes $f_i$ (resp.  $f_{i+1}$). As a result we have
  $||F_0 p_{\alpha'}||^2 - ||F_0 p_{\alpha}||^2 = f_{i+1} - f_i$ and the first step gives
  the desired upper bound.
\end{dem}

\begin{thm} \label{thm:AO}
  Let $S$ be a Woronowicz \Cst algebra and $p_1$ a central projection of $\hat S$ such
  that $Up_1U = p_1$ and $p_0 p_1 = 0$. Assume that the classical Cayley graph $\GG$ of
  $(S, p_1)$ is a directional tree and that $M_\gamma \neq 2$ for all $\gamma\in\Dir$. Let
  $F : H \to H\tens H$ be the isometry of Definition~\ref{df:isom_AO}. Then
  \begin{displaymath}
    F^* (\lambda\tens\rho)(x) F \equiv (\lambda,\rho)(x) ~\mathrm{mod}~ K(H) \text{,}
  \end{displaymath}
  for any $x \in S_\red\tensmax S_\red$. In particular $(S,\delta)$ has Property AO.
\end{thm}

\begin{dem}
  By symmetry one can assume that $x = a\tens 1$ with $a\in S_\red$ a coefficient of some
  $\gamma\in\Dir$. We put $G = (F_0^*F_0)^{-\frac12}$, so that $F = F_0 G$. We have then
  \begin{displaymath}
    aF^* - F^*(a\tens 1) = G [G^{-1}, a] F^* + G \Te [\Pe_{\ss\jok+}, a\tens 1] \text{.} 
  \end{displaymath}
  By the first statement of Lemma~\ref{lem:norm_AO}, the operator $G$ is compact. Hence it
  suffices to prove that $[G^{-1},a]$ and $\Te [\Pe_{\ss\jok+}, a\tens 1]$ are bounded.
  
  For the first commutator, we remark that $a p_\alpha = (p_{\alpha'} + p_{\alpha''}) a
  p_\alpha$ if $\gamma\tens\alpha = \alpha' \oplus \alpha''$. Moreover we have $Gp_\alpha =
  ||F_0p_\alpha||^{-1} p_\alpha$, so that 
  \begin{displaymath}
   [G^{-1},a] p_\alpha = (||F_0p_{\alpha'}|| - ||F_0p_{\alpha}||) p_{\alpha'} a p_\alpha
   + (||F_0p_{\alpha''}|| - ||F_0p_{\alpha}||) p_{\alpha''} a p_\alpha \text{.}
  \end{displaymath}
  Hence the result follows from the second statement of Lemma~\ref{lem:norm_AO}, after
  factoring out $||F_0p_{\alpha'}|| + ||F_0p_{\alpha}||$ from it. (In fact this even proves
  that $[G^{-1}, a]$ is compact.)
  
  For the second commutator, we will assume that $\dim\gamma > 1$: otherwise the proof is
  as easy as in the classical case. Denote by $(m_k)_k$ the sequence of quantum dimensions
  associated with $\gamma$. We use the Remarks~\ref{rque:ext_target}
  and~\ref{rque:reg_commut}.\ref{item:reg_asc_commut_ext} to write
  \begin{eqnarray*}
    ||p_\alpha \Te [\Pe_{\ss\jok+}, a\tens 1]||^2 \leq \makebox[-3cm]{} && \\[2ex]
    &\leq& \sum ||p_\alpha \Te (p_{\beta_1}\tens p_{\beta_2})||^2
    \times || [\Pe_{\ss\jok+}, a^*\tens 1] (p_{\beta_1}\tens p_{\beta_2})||^2 \\
    &\leq& C_a~ \sum \frac{m_{|\beta_1|} m_{|\beta_2|}} {m_{|\alpha|}} 
    \frac{m_{|\beta_2|-1}} {m_{|\alpha|-1} m_{|\beta_1|}} =
    \frac {C_a}{m_{|\alpha|} m_{|\alpha|-1}} \sum_{k=0}^{|\alpha|} m_k m_{k-1} \text{.}
  \end{eqnarray*}
  The last upper estimate is bounded because $m_k \sim \frac{a^{k+1}} {a-a^{-1}}$ for some $a >
  1$.
\end{dem}

\begin{rque}{Remark}
  Recall that $A_u(Q)$ is never amenable and that $A_o(Q)$ is amenable \iff $n=2$
  \cite{Banic:U(n)}. Hence the only case, up to free products, where Property AO is
  neither trivial nor proved by Theorem~\ref{thm:AO}, is the one of
  $A_u\big({1\atop0}{0\atop1}\big)$. Property AO may however be fulfilled in this case,
  too.
\end{rque}

\subsection{$KK$-theory}
\label{section:gamma}

The notion of $K$-amenability was first introduced by Cuntz \cite{Cuntz:kmoy} for discrete
groups: the aim was to give a simpler proof to a result of Pimsner and Voiculescu
\cite{PimsnerVoicu:freecross} calculating the $K$-theory of the reduced \Cst algebras of 
free groups. Cuntz proves that the $K$-theory of the reduced and full \Cst algebras of a
free group are the same, and gives in \cite{Cuntz:freeprod} a simple way to compute
it in the full case. 

Julg and Valette extended then the notion of $K$-amenability to the locally compact case
and established the $K$-amenability of locally compact groups acting on trees with
amenable stabilizers \cite{JulgValette:kmoy}. This includes the case of the free groups
acting on their Cayley graphs. To prove the $K$-amenability of a locally compact group $G$,
one has to construct an element $\alpha \in KK_G(\CC,\CC)$ using representations of $G$
that are weakly contained in the regular one, and then to prove that $\alpha$ is homotopic
to the unit element of $KK_G(\CC,\CC)$. In \cite{JulgValette:kmoy}, both of the steps are
carried out in a very geometric way. Moreover, it turns out that $\alpha$ can be
interpreted as the $\gamma$ element used to prove the Baum-Connes conjecture in this
context \cite{KaspSkand:buildings_novikov}.

\bigskip

We refer the reader to \cite{BaajSkand:kks,Vergnioux:these} for details about equivariant
$KK$-theory with respect to Hopf \Cst algebras, and we just recall the equivalent
characterizations of $K$-amenability for a discrete quantum group defined by its full and
reduced Woronowicz \Cst algebras $S$, $S_\red$:
\begin{enumerate} \setlength{\itemsep}{0ex}
  \renewcommand\theenumi{{\it\roman{enumi}}} \renewcommand\labelenumi{{\it\roman{enumi}}.}
\item \label{enum:kmoy_un} $\un \in KK_{\hat S}(\CC,\CC)$ can be represented by a triple
  $(E,\pi,F)$ such that the representation of $S$ on $E$ factors through $S_\red$.
\item \label{enum:kmoy_croise} For every \Cst algebra $A$ endowed with a coaction of $\hat
  S$, \newline $[\lambda_A] \in$ ${KK(A\rtimes S, A\rtimes_\red S)}$ is invertible.
\item \label{enum:kmoy_lambda} $[\lambda] \in KK(S, S_\red)$ is invertible.
\item \label{enum:kmoy_alpha} There exists $\alpha \in KK(S_\red,\CC)$ such that
  $\lambda^*(\alpha) = [\varepsilon] \in KK(S,\CC)$.
\end{enumerate}
In this subsection we explain how to construct an element $\alpha \in KK(S_\red,\CC)$ from
the quantum Cayley graph of a free quantum group. It is the natural quantum generalization
of the Julg-Valette element mentioned above. It has index $1$, however further work is
needed to determine whether $\lambda^*(\alpha) = [\varepsilon]$.

\begin{thm}
  Let $S$ be a Woronowicz \Cst algebra and $p_1$ a central projection of $\hat S$ such
  that $Up_1U = p_1$ and $p_0 p_1 = 0$. Assume that the classical Cayley graph $\GG$ of
  $(S, p_1)$ is a directional tree and that $M_\gamma \neq 2$ for all $\gamma\in\Dir$.
  Then $\T p_{\ss++} : K_g \to H$ and $\T(Rs)^* : K_\infty \to H$ commute to the actions
  of $S_\red$ modulo compact operators. In particular $\T (p_{\ss++} + (Rs)^*)$ defines an
  element $\alpha \in KK(S_\red,\CC)$ of index $1$.
\end{thm}

\begin{dem}
  In this proof we will denote by $p_{\geq k_0}$ the sum of the projections $p_k$ with
  $k\geq k_0$. We have $\T p_{\ss++} = \T p_{\ss\jok+}$ and the target operator $\T$
  intertwines the actions of $S_\red$, hence it is enough to prove that $p_{\ss\jok+}$
  commutes to $S_\red\tens 1$ up to compact operators. If $a \in S_\red$ is a coefficient
  of some $\gamma\in\Dir$, this results directly from Lemma~\ref{lem:reg_asc_commut}:
  because $(m_k^{-1})$ is decreasing and $p_{k'}ap_k$ vanishes as soon as $|k-k'| \neq 1$,
  we have
  \begin{displaymath}
    || [p_{\ss\jok+}, a\tens 1] (p_k\tens\id) || \leq C_a m_k^{-1} ~\Longrightarrow~ 
    || [p_{\ss\jok+}, a\tens 1] (p_{\geq k_0}\tens\id) || \leq 2 C_a m_{k_0}^{-1} \text{.}
  \end{displaymath}
  This proves that $[p_{\ss\jok+}, a\tens 1]$ is compact, and the general result follows
  because the coefficients $a$ of the corepresentations $\gamma\in\Dir$ span the \Cst
  algebra $S_\red$.
  
  For the case of $(Rs)^*$ we will use the proof of Theorem~\ref{thm:infinite_reg}. Thanks
  to the hypothesis we can take into account the simplification of
  Remark~\ref{rque:simplif_infinite_reg}: we avoid the use of the projections $Q_{k_0,l}$
  by taking for $\rho_k$ the remainder of $(\sum m_k^{-1})$ instead of $(\sum m_k^{-2})$.
  Equation~(\ref{eq:infinite_reg_conv}) reads then
  \begin{displaymath}
    || (\pi_\infty(a)R - Rp_{\ss+-}(a\tens 1)) (p_k\tens\id)|| \leq (2p+1) C_a \rho_k
    \text{.} 
  \end{displaymath}
  We notice also that $\rho_k$ is again equivalent to a multiple of $m_k^{-1}$ because
  $(m_k)$ grows geometrically. To conclude we use Lemma~\ref{lem:reg_asc_commut} and the
  fact that $(a\tens 1)$ commutes to $\Theta$: up to a change of the constant $C_a$, we
  obtain $|| (\pi_\infty(a)Rs - Rs(a\tens 1)) (p_k\tens\id)|| \leq C_a m_k^{-1}$. Summing
  over $k\geq k_0$ we obtain an inequality showing that $(\pi_\infty(a)Rs - Rs(a\tens 1))$
  is compact:
  \begin{displaymath}
    || (\pi_\infty(a) Rs - Rs(a\tens 1)) (p_{\geq k_0}\tens \id)|| \leq C_a \rho_{k_0}
    \text{.} 
  \end{displaymath}
  
  Finally $\T (p_{\ss++} + (Rs)^*)$ defines an element $\alpha \in KK(S_\red,\CC)$ of
  index $1$ because $\T : K_{\ss++} \to (1-p_0)H$ is invertible by
  Proposition~\ref{prp:but_inj}, as well as $p_{\ss++} + (Rs)^* : K_g \oplus K_\infty \to
  K_{\ss++}$ by Theorem~\ref{thm:proj}, Proposition~\ref{prp:RS} and
  Theorem~\ref{thm:proj2}.
\end{dem}


\small
\begin{thebibliography}{10}

\bibitem{AkeOstrand:tensorfree}
Charles~A. Akemann and Phillip~A. Ostrand.
\newblock On a tensor product {$C\sp*$}-algebra associated with the free group
  on two generators.
\newblock {\em J. Math. Soc. Japan}, 27(4):589--599, 1975.

\bibitem{BaajSkand:kks}
Saad Baaj and Georges Skandalis.
\newblock {$C\sp \ast$}-alg\`ebres de {H}opf et th\'eorie de {K}asparov
  \'equivariante.
\newblock {\em $K$-Theory}, 2(6):683--721, 1989.

\bibitem{BaajSkand:unit}
Saad Baaj and Georges Skandalis.
\newblock Unitaires multiplicatifs et dualit\'e pour les produits crois\'es de
  {$C\sp *$}-alg\`ebres.
\newblock {\em Ann. Sci. \'Ecole Norm. Sup. (4)}, 26(4):425--488, 1993.

\bibitem{Banic:O(n)_cras}
Teodor Banica.
\newblock Th\'eorie des repr\'esentations du groupe quantique compact libre
  {${\rm O}(n)$}.
\newblock {\em C. R. Acad. Sci. Paris S\'er. I Math.}, 322(3):241--244, 1996.

\bibitem{Banic:U(n)}
Teodor Banica.
\newblock Le groupe quantique compact libre {${\rm U}(n)$}.
\newblock {\em Comm. Math. Phys.}, 190(1):143--172, 1997.

\bibitem{Cuntz:freeprod}
Joachim Cuntz.
\newblock The {$K$}-groups for free products of {$C\sp{\ast} $}-algebras.
\newblock In {\em Operator algebras and applications, Part I (Kingston, Ont.,
  1980)}, volume~38 of {\em Proc. Sympos. Pure Math.}, pages 81--84. Amer.
  Math. Soc., Providence, R.I., 1982.

\bibitem{Cuntz:kmoy}
Joachim Cuntz.
\newblock {$K$}-theoretic amenability for discrete groups.
\newblock {\em J. Reine Angew. Math.}, 344:180--195, 1983.

\bibitem{JulgValette:kmoy}
Pierre Julg and Alain Valette.
\newblock {$K$}-theoretic amenability for {${\rm SL}\sb{2}({\mathbb
  Q}\sb{p})$}, and the action on the associated tree.
\newblock {\em J. Funct. Anal.}, 58(2):194--215, 1984.

\bibitem{KaspSkand:buildings_novikov}
Gennadi~G. Kasparov and Georges Skandalis.
\newblock Groups acting on buildings, operator {$K$}-theory, and {N}ovikov's
  conjecture.
\newblock {\em $K$-Theory}, 4(4):303--337, 1991.

\bibitem{KustVaes:lcqg}
Johan Kustermans and Stefaan Vaes.
\newblock Locally compact quantum groups.
\newblock {\em Ann. Sci. \'Ecole Norm. Sup. (4)}, 33(6):837--934, 2000.

\bibitem{nikolskij:toeplitz}
Nikolai~K. Nikol'ski\u\i.
\newblock {\em Hardy, Hankel and Toeplitz}, volume~1 of {\em Operators,
  functions and systems}.
\newblock AMS Mathematical surveys and monographs, 1992.

\bibitem{Ozawa:solid}
Narutaka Ozawa.
\newblock Solid von {Neumann} algebras.
\newblock Preprint, 2003.

\bibitem{PimsnerVoicu:freecross}
Mihai~V. Pimsner and Dan~V. Voiculescu.
\newblock {$K$}-groups of reduced crossed products by free groups.
\newblock {\em J. Operator Theory}, 8(1):131--156, 1982.

\bibitem{Serre:arbres}
Jean-Pierre Serre.
\newblock {\em Arbres, amalgames, {${\rm SL}\sb{2}$}}.
\newblock Soci\'et\'e Math\'ematique de France, Paris, 1977.
\newblock Avec un sommaire anglais, R\'edig\'e avec la collaboration de Hyman
  Bass, Ast\'erisque, No. 46.

\bibitem{Skand:knucl}
Georges Skandalis.
\newblock Une notion de nucl\'earit\'e en {$K$}-th\'eorie (d'apr\`es {J}.\
  {C}untz).
\newblock {\em $K$-Theory}, 1(6):549--573, 1988.

\bibitem{DaeleWang:univ}
Alfons Van~Daele and Shuzhou Wang.
\newblock Universal quantum groups.
\newblock {\em Internat. J. Math.}, 7(2):255--263, 1996.

\bibitem{Vergnioux:these}
Roland Vergnioux.
\newblock {\em $KK$-th\'eorie \'equivariante et op\'erateur de Julg-Valette
  pour les groupes quantiques}.
\newblock Th\`ese de doctorat, Universit\'e Paris 7, 2002.

\bibitem{Vergnioux:free}
Roland Vergnioux.
\newblock ${K}$-amenability for amalgamated free products of amenable discrete
  quantum groups.
\newblock To appear in {\em Journal of Functional Analysis}, 2003.

\bibitem{Wang:freeprod}
Shuzhou Wang.
\newblock Free products of compact quantum groups.
\newblock {\em Comm. Math. Phys.}, 167(3):671--692, 1995.

\bibitem{Woro:cmp}
Stanis\l aw~L. Woronowicz.
\newblock Compact matrix pseudogroups.
\newblock {\em Comm. Math. Phys.}, 111(4):613--665, 1987.

\bibitem{Woro:dual}
Stanis\l aw~L. Woronowicz.
\newblock Tannaka-{K}re\u\i n duality for compact matrix pseudogroups.
  {T}wisted {${\rm SU}(N)$} groups.
\newblock {\em Invent. Math.}, 93(1):35--76, 1988.

\bibitem{Woro:houches}
Stanis\l aw~L. Woronowicz.
\newblock Compact quantum groups.
\newblock In {\em Sym\'etries quantiques (Les Houches, 1995)}, pages 845--884.
  North-Holland, Amsterdam, 1998.

\end{thebibliography}

\normalsize

\bigskip

\begin{center}
  \parbox{10cm}{\noindent
    Mathematisches Institut, Westf\"alische Wilhelms-Universit\"at \\
    Einsteinstr. 62, D--48149 M\"unster \\
    e-mail: {\tt vergniou@math.uni-muenster.de} \\
    3. Februar 2004}
\end{center}

\end{document}